\documentclass[preprint]{elsarticle}
\usepackage{mathrsfs}
\usepackage{subcaption}
\usepackage{bm}\makeatletter \oddsidemargin -.1in \evensidemargin -.1in
\newcommand{\singlespacing}{\let\CS=\@currsize\renewcommand{\baselinestretch}{1}\tiny\CS}
\newcommand{\doublespacing}{\let\CS=\@currsize\renewcommand{\baselinestretch}{1.27}\tiny\CS}
\makeatletter
\newcommand{\setword}[2]{%
	\phantomsection
	#1\def\@currentlabel{\unexpanded{#1}}\label{#2}%
}
\makeatother

\newtheorem{definition}{Definition}[section]

\newtheorem{condition}{Condition}[section]
\newtheorem{example}{Example}[section]
\newtheorem{theorem}{Theorem}[section]

\newproof{proof}{Proof}
\newtheorem{lemma}{Lemma}[section]

\newtheorem{remark}{Remark}[section]

\newcommand{\osc}{\operatorname{osc}}
\newcommand{\dimB}{\dim_{B}}
\newcommand{\overdimB}{\overline{\dim}_{B}}
\newcommand{\underdimB}{\underline{\dim}_{B}}
\usepackage{xcolor}
\usepackage{microtype}
\usepackage{times}
\usepackage{graphicx,epsfig}
\usepackage{amsmath}
\usepackage{amssymb}
\usepackage{enumitem}
\usepackage{mathtools}
\usepackage{float}
\everymath{\displaystyle}
\usepackage[linkcolor=blue, urlcolor=blue, citecolor=blue,
colorlinks, bookmarks]{hyperref}
\date{}

\numberwithin{equation}{section}
\journal{Journal...}

\begin{document}
	
	\begin{frontmatter}
		
		\title {\textbf{\Large Partial Self-Similar Fractal Interpolation Functions and Surfaces}}

		\author[1]{Pavjeet Singh}
		\ead{pavjeetsingh23@gmail.com}
		
		\author[2]{Megha Pandey\corref{mycorrespondingauthor}}
		\cortext[mycorrespondingauthor]{Corresponding author}
		\ead{meghapandey1071996@gmail.com}

		\author[1]{S. K. Katiyar}
		\ead{katiyarsk@nitj.ac.in}
		\author[3]{Tanmoy Som}
		\ead{tsom.apm@iitbhu.ac.in}
		\address[1]{DR BR Ambedkar National Institute of Technology, Jalandhar, Punjab, 144008, India}
		\address[2]{School of Mathematics, Northwest University Xi'an, Shaanxi Province, 710127, PR China}
		\address[3]{Department of Mathematical Sciences, Indian Institute of Technology (Banaras Hindu University) Varanasi, UP, 221005, India}
		
		\begin{abstract}
			In this paper, we study a class of fractal interpolation functions with partial self-similarity. We first study fractal interpolation surfaces defined on rectangular domains, analyze the stability and sensitivity of the proposed construction to parameter perturbations. In addition, we investigate the partial $\alpha$-fractal function in both the univariate and bivariate cases, as well as the invariant measure in the univariate setting whose support is the graph of the fractal interpolation function. Further, under suitable additional conditions, we derive bounds for the box-counting dimension of the graph of the univariate PFIF. The results presented in this paper are backed up by illustrative examples.
			
			
			
		\end{abstract}

	\end{frontmatter}
	
	\section{Introduction}
	Fractal interpolation functions (FIFs), introduced by Barnsley \cite{barnsley1986fractal,barnsley2014fractals}, provide a useful framework for constructing interpolants that can capture irregular and complicated structures present in the data. In contrast to the classical interpolation methods, FIFs are capable of modeling highly irregular and oscillatory phenomena, making them useful in several areas such as approximation theory, signal processing, and image compression \cite{mazel1994representation, david1992using}. Following the pioneering work of Barnsley, the theory of FIFs has developed rapidly over the past few decades. Various generalizations have been introduced, including hidden variables viz FIFs \cite{barnsley1989hidden,chand2003hidden}, recurrent FIFs \cite{barnsley1989recurrent}, bilinear FIFs \cite{barnsley2015bilinear}, spline and Hermite FIFs \cite{chand2006generalized}, as well as nonlinear fractal interpolation function \cite{ri2018new,ri2019new}, rational cubic spline FIFs and their constrained and shape-preserving aspects \cite{katiyar2019new}. Interpolative operators linking multivalued fractal structure \cite{prithvi2022interpolative}, while nonconventional iterated function systems (IFSs) have been used to further extend the framework of fractal interpolation \cite{prithvi2023revisiting}.
	
	The extension of fractal interpolation theory to higher dimensions has received considerable attention in recent years. Several authors have studied bivariate/multivariate fractal interpolation functions and associated fractal operators on rectangular grids \cite{bouboulis2007fractal,ruan2015fractal,pandey2026generalized,agrawal2022bivariate}. These developments provide flexible tools for modeling irregular surfaces arising in applied sciences and engineering. Another important development in fractal approximation theory is the study of $\alpha$-fractal functions and fractal operators. In perturbation techniques, $\alpha$-fractal functions preserve interpolation properties while incorporating fractal characteristics \cite{navascues2005fractal}. Several approximation-theoretic and geometric properties of such functions have been investigated in \cite{pandey2022fractal, pandey2024set}. Analytical properties of fractal interpolation functions such as smoothness, differentiability, integration, and shape-preserving properties have also been widely investigated in \cite{barnsley1989calculus,viswanathan2014fractal}.  In particular, box-counting dimension estimates for FIFs and related fractal functions have been studied in \cite{barnsley2015bilinear, falconer2013fractal, akhtar2017box}. 
	
	The graph of FIFs constructed in \cite{barnsley1986fractal} exhibits self-similarity in nature. However, self-similarity in a real object need not be global. For example, consider a mountainous landscape. One part of the landscape may contain a
	pattern of ridges and valleys similar to another particular part, while it may look quite different from the landscape as a whole. In such a situation, it is natural to allow different parts of the interpolation function to be generated from different subintervals rather than the whole
	interpolation interval. This idea was introduced by Luor~\cite{luor2018fractal} through fractal interpolation functions with partial self-similarity. In this setting, the $k$th interpolation interval $I_k$ is associated with a subinterval $J_k$, and the corresponding map has the form
	\[J_k \longrightarrow I_k,\]
	where $J_k$ may depend on $k$. Thus, different pieces of the interpolation functions can be generated from different parts of its graph.
	
	Recurrent fractal interpolation functions are closely related to fractal interpolation functions with partial self-similarity. In the recurrent setting, each interpolation interval $I_i$ is associated with an interval $D_i$, and the graph over $I_i$ is generated from those component graphs corresponding to the basic intervals contained in $D_i$. Similarly, in Luor's partial self-similar construction, each interval $I_k$ is associated with an interval $J_k$, and the graph over $J_k$ is used to generate the graph over $I_k$. Consequently, when $J_k$ is decomposed into the basic interpolation intervals, the partial self-similar construction also gives rise to a graph-directed structure of recurrent type.
	The distinction is therefore mainly in the formulation with emphasis of the two approaches. The recurrent framework describes the dependence through the associated connection structure, whereas Luor's formulation emphasizes the local self-similar relation between the graph over $J_k$ and the graph over $I_k$. An important feature of Luor's construction is that $J_k$ is not required to be larger than $I_k$ and the map $L_k:J_k\to I_k$ need not be contractive. In the present work, we extend this partial self-similar viewpoint to the bivariate setting, where analogous restrictions on the horizontal mappings are not imposed.
	
	Further, the use of different subintervals also leads to a structure different from the usual IFS. In an ordinary IFS, all the maps act on the same space, and the invariant measure is defined through these maps on that common space. For a fractal interpolation function with partial self-similarity, however, different interpolation intervals may be generated from different parts of the graph. Hence, the relations between the different parts of the graph are naturally described by a graph-directed system. In this system, an edge
	\[G_r \longrightarrow G_k\]
	indicates that the graph over $I_r$ can be used to construct the graph over $I_k$. This graph-directed structure naturally leads to a different invariant
	measure problem. Rather than using a single invariant measure associated with maps acting on the whole graph, we consider a measure on each part of the graph, with these measures related through the graph-directed transitions. We develop this invariant-measure construction for fractal interpolation functions with partial self-similarity in the univariate setting. We also study the box dimension of the corresponding interpolation graph, which describes its geometric complexity.
	
	The box dimension of recurrent fractal interpolation functions was initially studied by Barnsley et al. \cite{barnsley1989recurrent}. Later, Ruan et al. \cite{ruan2021existence} generalized this framework by removing the requirement that the horizontal maps $L_i:D_i\to I_i$ be contractive and obtained the box dimension for self-affine RFIFs in terms of the spectral radius of an associated weighted connection matrix. Under the additional conditions considered here, the univariate PFIF leads to the same dimension characterization, although the proof is developed directly within the PFIF framework.
	
	Motivated by these considerations, we extend partial self-similarity to the bivariate setting. In particular, instead of considering only intervals, we consider rectangles and allow the construction
	\[R_k\times Q_m \longrightarrow I_k\times J_m,\]
	where the choice of $R_k$ and $Q_m$ may depend on the corresponding interpolation intervals. This provides a framework for studying fractal interpolation surfaces in which the partial self-similar structure can vary from one interpolation interval to another, while the invariant measure and box dimension results are developed in the univariate setting.
	
	\subsection{Delineation}
	The paper is organized as follows. In the next section, we review the necessary preliminaries and recall Luor's construction, which is used in our paper. Section~\ref{sec3} is devoted to the extension to bivariate fractal interpolation surfaces, along with illustrative examples. In Section~\ref{sec4}, we study partial $\alpha$-fractal interpolation functions associated with the fractal interpolation functions and derive their error estimates. Section~\ref{sec5} discusses the stability and sensitivity analysis of the fractal interpolation functions. In Section~\ref{sec6}, we investigate the invariant measure associated with the system. Further, Section~\ref{sec7} presents the box-counting dimension analysis. Finally, the paper is concluded in Section~\ref{sec8} with some future directions.

	\section{Notations and Preliminaries}
	\subsection{Notations}
	Throughout this paper, we shall use the following notations:
	\begin{itemize}
		\item $\mathbb{R}$: set of real numbers.
		\item $\mathbb{N}$: set of natural numbers.
		\item $I, J$: unit intervals $[0,1]$.
		\item $I\times J$: Cartesian product of $I \text{ and } J$.
		\item $\lvert I \rvert$: length of the interval $I$.
		\item $\Sigma_K=\{1,2,\ldots,K\},~~ \Sigma_{K,0}=\{0,1,\ldots, K\},~~ \partial \Sigma_{K,0}=\{0,K\}, \text{ and } \text{int} \Sigma_{K,0}=\{1,2,\ldots,K-1\}$, where $K\in \mathbb{N}$.
	\end{itemize}

	\subsection{Preliminaries}
	Let $(\boldsymbol{Y}, d_{\boldsymbol{Y}})$ be a complete metric space and $I=[0,1] $ and $J=[0,1]$ be two unit intervals. Taking the interpolation data as 
	\[\Delta =\left\{(t_i, s_j, y_{ij})\in I\times J\times \boldsymbol{Y} : i=0,1,\ldots, N; j=0,1,\ldots, M \right\},\] 
	with $0=t_0<\cdots < t_N=1$ and $0=s_0<\cdots< s_M=1$ and $N,M \geq 2$. \\
	Let $\boldsymbol{X}=I\times J\times \boldsymbol{Y}$ and $\theta >0$. Define a metric $d_{\boldsymbol{X}}$ on $\boldsymbol{X}$ by 
	\[d_{\boldsymbol{X}}((t,s,y),(t^*,s^*,y^*))=\lvert t-t^*\rvert + \lvert s-s^*\rvert +\theta d_{\boldsymbol{Y}}(y,y^*) \text{ for all } (t,s,y), (t^*,s^*,y^*)\in \boldsymbol{X}.\]
	Then, $(\boldsymbol{X}, d_{\boldsymbol{X}})$ is a complete metric space.\\

	\noindent
	Throughout this paper, we denote by $\mathcal{C}(I\times J, \boldsymbol{Y})$ the set of all continuous functions defined on $I\times J$ and $\mathcal{C}_{\Delta}(I\times J, \boldsymbol{Y})$ be the collection of all those continuous functions which interpolate the points in $\Delta$. The graph of $f\in \mathcal{C}(I\times J, \boldsymbol{Y})$ is defined by 
	\[Gr(f)=\{(t,s,f(t,s)): t\in I \text{ and } s\in J\}.\]
	Let $\mathbf{G}(I\times J, \boldsymbol{Y})=\left\{Gr(f): f\in \mathcal{C}(I\times J, \boldsymbol{Y})\right\}$ and $\mathbf{G}_{\Delta}=\left\{Gr(f): f\in \mathcal{C}_{\Delta}(I\times J, \boldsymbol{Y})\right\}$. For $A\in \mathbf{G}(I\times J, \boldsymbol{Y})$, 
	denote $f_{[A]}\in \mathcal{C}(I\times J, \boldsymbol{Y})$, a unique function such that $A=Gr\left(f_{[A]}\right)$. Consider the sup metric $d_{\infty}$ on $\mathcal{C}(I\times J, \boldsymbol{Y})$,
	\[d_{\infty}(f,g)=\max_{(t,s)\in I\times J}d_{\boldsymbol{Y}}\left(f(t,s), g(t,s)\right) \text{ for } f,g\in \mathcal{C}(I\times J, \boldsymbol{Y}).\]
	It is proven that $\left(\mathcal{C}(I\times J, \boldsymbol{Y}), d_{\infty}\right)$ is a complete metric space and since $\left(\mathcal{C}_{\Delta}(I\times J, \boldsymbol{Y}), d_{\infty}\right)$ is a closed subspace of $\left(\mathcal{C}(I\times J, \boldsymbol{Y}), d_{\infty}\right)$, therefore it is also a complete metric space.
	For $A,B\in \mathbf{G}(I\times J, \boldsymbol{Y})$, define
	\[d_{\mathbf{G}}(A,B)=d_{\infty}\left(f_{[A]}, f_{[B]}\right).\]
	Then, note that $\left(\mathbf{G}(I\times J, \boldsymbol{Y}), d_{\mathbf{G}}\right)$ and $\left(\mathbf{G}_{\Delta}(I\times J, \boldsymbol{Y}), d_{\mathbf{G}}\right)$ both are complete metric spaces.
	
	\begin{definition}[see \cite{falconer2013fractal}]
		Let $G$ be any non-empty bounded subset of $\mathbb{R}^n$ and let $N_\delta(G)$ be the minimum number of square boxes of side $\delta$ parallel to the axes which can cover $G$. The upper and lower box-counting dimensions of $G$ are respectively defined as 
		\[\overdimB G	=\limsup_{\delta \to 0}\frac{\log N_\delta(G)}{-\log\delta},\]
		and 
		\[\underdimB G	=\liminf_{\delta \to 0}\frac{\log N_\delta(G)}{-\log\delta} .\]
		If both limits exist and are equal, then the box dimension of $G$ is $\dim_B G$ and is defined as
		\[\dimB G	=\lim_{\delta \to 0}\frac{\log N_\delta(G)}{-\log\delta} .\]
	\end{definition}

	\section{Construction of bivariate fractal interpolation functions}\label{sec3}
	
	\noindent
	For each $I_{k}=[t_{k-1},t_k]$, define $R_{k}=\left[t_{j(k)},t_{l(k)}\right]$ such that $j(k)<l(k)$ and for each $J_m=[s_{m-1},s_m]$, define $Q_m=\left[s_{w(m)}, s_{z(m)}\right]$ such that $w(m)<z(m)$, where $k\in \Sigma_{N}, ~ m\in \Sigma_{M} \text{ and } j,l :\Sigma_{N}\rightarrow \Sigma_{N,0},~ w,z:\Sigma_{M}\rightarrow \Sigma_{M,0}$.
	
	\noindent
	Next we introduce the functions $U_k,~ V_m,$ and $P_{km}$ for $k=1,\ldots, N$ and $m=1,\ldots, M$ as follows:
	\begin{enumerate}[label=(\Roman*)]
		\item $U_k:R_k\rightarrow I_k$ is a homeomorphism such that $U_{k}\left(t_{j(k)}\right)=t_{k-1}$ and $U_k\left(t_{l(k)}\right)=t_k.$
		\item  $V_m: Q_m\rightarrow J_m$ is also a homeomorphism such that $V_m\left(s_{w(m)}\right)=s_{m-1}$ and $V_m \left(s_{z(m)}\right)=s_m$.
		\item $P_{km}: R_k\times Q_m\times \boldsymbol{Y}\rightarrow \boldsymbol{Y}$ is a continuous function such that 
		\begin{align}\label{cond}
			\begin{cases}
				P_{km}\left(t_{j(k)}, s_{w(m)}, y_{j(k)w(m)}\right)= y_{(k-1)(m-1)},  \quad P_{km}\left(t_{j(k)}, s_{z(m)}, y_{j(k)z(m)}\right)= y_{(k-1)m},\\
				P_{km}\left(t_{l(k)}, s_{w(m)}, y_{l(k)w(m)}\right)= y_{k(m-1)}, \quad \hspace{0.6cm}P_{km}\left(t_{l(k)}, s_{z(m)}, y_{l(k)z(m)}\right)=y_{km}.
			\end{cases}
		\end{align}
	\end{enumerate}
	For each $k$ and $m$, assume that the condition
	\begin{equation}\label{contractive}
		d_{\boldsymbol{Y}}\left(P_{km}(t,s,y), P_{km}(t,s,y^{*})\right)\leq \eta_{km}d_{\boldsymbol{Y}}(y,y^{*}) \text{ for all } t\in R_k,~ s\in Q_m, \text{ and } y,y^*\in \boldsymbol{Y},
	\end{equation}
	holds for some $0<\eta_{km}<1$. In addition, assume that $\{P_{km}\}$ satisfies the
	following join-up conditions. For every
	$f\in\mathcal{C}_{\Delta}(I\times J,\boldsymbol{Y})$, we have
	\begin{align}
		P_{km}\left(
		t_{l(k)},v,f(t_{l(k)},v)
		\right) =
		P_{(k+1)m}\left(
		t_{j(k+1)},v,f(t_{j(k+1)},v)
		\right),
	\end{align}
	for every $k=1,\ldots,N-1$, $m=1,\ldots,M$, and $v\in Q_m$.
	Similarly,
	\begin{align*}
		P_{km}\left(
		u,s_{z(m)},f(u,s_{z(m)})
		\right)  =
		P_{k(m+1)}\left(
		u,s_{w(m+1)},f(u,s_{w(m+1)})
		\right),
	\end{align*}
	for every $k=1,\ldots,N$, $m=1,\ldots,M-1$, and $u\in R_k$.  Define $W_{km}:R_k\times Q_m\times \boldsymbol{Y}\rightarrow I_{k}\times J_m\times \boldsymbol{Y}$ by 
	\begin{equation}\label{Wkm}
		W_{km}(t,s,y)=\left(U_{k}(t), V_{m}(s), P_{km}(t,s,y)\right) \text{ for all } t\in R_{k},~ s\in Q_m, \text{ and } y\in \boldsymbol{Y}.
	\end{equation}
	Observe that
	\begin{align*}
		W_{km}\left(t_{j(k)}, s_{w(m)}, y_{j(k)w(m)}\right)=& \left(t_{k-1}, s_{m-1}, y_{(k-1)(m-1)}\right) &\quad W_{km}\left(t_{j(k)}, s_{z(m)}, y_{j(k)z(m)}\right)=& \left(t_{k-1}, s_{m}, y_{(k-1)m}\right)\\
		W_{km}\left(t_{l(k)}, s_{w(m)}, y_{l(k)w(m)}\right)=& \left(t_{k}, s_{m-1}, y_{k(m-1)}\right) &\quad W_{km}\left(t_{l(k)}, s_{z(m)}, y_{l(k)z(m)}\right)=& \left(t_{k}, s_{m}, y_{km}\right).
	\end{align*}
	For $A\in \mathbf{G}_{\Delta}(I\times J, \boldsymbol{Y})$ and $k=1,\ldots, N$ and $m=1,\ldots, M,$ define $A_{R_kQ_m}=A\cap (R_k\times Q_m\times \boldsymbol{Y})$, then 
	\begin{align}
		& W_{km}\left(A_{R_kQ_m}\right)=\left\{\left(U_k(t), V_m(s), P_{km}\left(t,s,f_{[A]}(t,s)\right)\right): t\in R_k \text{ and } s\in Q_m\right\}, \label{eq1}\\
		\implies& W_{km}\left(A_{R_kQ_m}\right)=\left\{\left(t,s,P_{km}\left(U_{k}^{-1}(t), V_{m}^{-1}(s),f_{[A]}\left(U_{k}^{-1}(t), V_{m}^{-1}(s)\right)\right)\right): t\in I_k \text{ and } s\in J_m\right\}. \label{eq2}
	\end{align}
	Define a function $h_{km}:I_{k}\times J_{m}\rightarrow \boldsymbol{Y}$ by
	\begin{equation}\label{GWkm}
		h_{km}(t,s)=P_{km}\left(U_{k}^{-1}(t), V_{m}^{-1}(s),f_{[A]}\left(U_{k}^{-1}(t), V_{m}^{-1}(s)\right)\right) \text{ for } t\in I_k \text{ and } s\in J_m.
	\end{equation}
	Then, it is easy to observe that $h_{km}$ is a continuous function on $I_k\times J_m$ and $Gr(h_{km})=W_{km}\left(A_{R_kQ_m}\right)$.\\
	Since $U^{-1}_{k}(t_{k-1})=t_{j(k)}, ~ U^{-1}_{k}(t_k)=t_{l(k)},~ V_{m}^{-1}(s_{m-1})=s_{w(m)}, ~V_{m}^{-1}(s_m)=s_{z(m)}$ and $f_{[A]}\in \mathcal{C}_{\Delta}(I\times J, \boldsymbol{Y}),$ we have 
	\begin{align*}
		h_{km}(t_{k-1},s_{m-1})=& y_{(k-1)(m-1)}, &\quad h_{km}(t_{k-1},s_{m})=&y_{(k-1)m},\\
		h_{km}(t_k, s_{m-1})=& y_{k(m-1)}, &\quad h_{km}(t_k,s_m)=& y_{km}.
	\end{align*}
	
	\noindent
	Define
	\begin{equation}\label{omega}
		\boldsymbol{\Omega}(A)=\bigcup_{k=1}^{N}\bigcup_{m=1}^{M}W_{km}\left(A_{R_kQ_m}\right).
	\end{equation}
	We have the following lemma.
	
	\begin{lemma}\label{Lemma1}
		For $A\in \mathbf{G}_{\Delta}(I\times J, \boldsymbol{Y})$, we get $\boldsymbol{\Omega}(A)\in \mathbf{G}_{\Delta}(I\times J, \boldsymbol{Y})$ and 
		\begin{equation}\label{lem1}
			f_{[\boldsymbol{\Omega}(A)]}(t,s)=P_{km}\left(U_{k}^{-1}(t), V_{m}^{-1}(s), f_{[A]}\left(U_{k}^{-1}(t),V_{m}^{-1}(s)\right)\right),
		\end{equation}
		for $ t\in I_k, ~ s\in J_m, \text{ where } k\in \Sigma_{N} \text{ and } m \in \Sigma_{M}.$
	\end{lemma}
	
	\begin{proof}
		We first verify that the functions agree on the common vertical and horizontal edges. Let $v=V_m^{-1}(s)\in Q_m$ and $s\in J_m$. Since $U_k^{-1}(t_k)=t_{l(k)} \text{ and } U_{k+1}^{-1}(t_k)=t_{j(k+1)}$,
		we obtain
		\begin{align*}
			h_{km}(t_k,s)&=P_{km}\left(t_{l(k)},v,f_{[A]}(t_{l(k)},v)\right)\\
			&=P_{(k+1)m}\left(t_{j(k+1)},v,f_{[A]}(t_{j(k+1)},v)\right)\\
			&=h_{(k+1)m}(t_k,s).
		\end{align*}
		Therefore, $h_{km}(t_k,s)=h_{(k+1)m}(t_k,s) \text{ for every }s\in J_m$. Similarly, let $u=U_k^{-1}(t)\in R_k$. Since $V_m^{-1}(s_m)=s_{z(m)}\text{ and } V_{m+1}^{-1}(s_m)=s_{w(m+1)}$, we have
		\begin{align*}
			h_{km}(t,s_m)&=P_{km}\left(u,s_{z(m)},f_{[A]}(u,s_{z(m)})\right)\\
			&=P_{k(m+1)}\left(u,s_{w(m+1)},f_{[A]}(u,s_{w(m+1)})\right)\\
			&=h_{k(m+1)}(t,s_m).
		\end{align*}
		Hence, $h_{km}(t,s_m)=h_{k(m+1)}(t,s_m) \text{ for every }t\in I_k$. Define a map $h: I\times J\rightarrow \boldsymbol{Y}$ such that 
		\[h(t,s)=h_{km}(t,s),\qquad (t,s)\in I_k\times J_m.\]
		Since the rectangles $I_k\times J_m$ form a
		finite closed cover of $I\times J$, and $h_{km}$ is continuous on
		each rectangle, the pasting lemma implies that $h$ is continuous
		on $I\times J$.
		Now we prove that $h(t_k,s_m)=y_{km}$ for all $(t_k,s_m,y_{km})\in \Delta$. For this take $(t_k,s_m,y_{km})\in \Delta$ arbitrary, then four cases will arise:\\
		\textbf{Case 1:} If $t_k\in I_{k}$ and $s_m\in J_{m}$, then we have
		\[h(t_k,s_m)=h_{km}(t_k,s_m)=y_{km}.\]
		\textbf{Case 2:} If $t_k\in I_{k+1}$ and $s_m\in J_{m}$, then we have
		\[h(t_k,s_m)=h_{(k+1)m}(t_k,s_m)=y_{km}.\]
		\textbf{Case 3:} If $t_k\in I_{k}$ and $s_m\in J_{m+1}$, then we have
		\[h(t_k,s_m)=h_{k(m+1)}(t_k,s_m)=y_{km}.\]
		\textbf{Case 4:} If $t_k\in I_{k+1}$ and $s_m\in J_{m+1}$, then we have
		\[h(t_k,s_m)=h_{(k+1)(m+1)}(t_k,s_m)=y_{km}.\]
		Thus, $h(t_k,s_m)=y_{km}$ for all $k \in \Sigma_{N}$ and $m\in \Sigma_{M}$. Hence, $h\in \mathcal{C}_{\Delta}(I\times J, \boldsymbol{Y})$.
		Next we show that $Gr(h)=\boldsymbol{\Omega}(A)$. For this, we have
		\begin{align*}
			Gr(h)=&\left\{(t,s,h(t,s)): t\in I \text{ and } s\in J\right\}\\
			=&\bigcup_{k=1}^{N}\bigcup_{m=1}^{M}\left\{(t,s,h(t,s)): t\in I_k \text{ and } s\in J_m \right\}\\
			=& \bigcup_{k=1}^{N}\bigcup_{m=1}^{M} \left\{(t,s,h_{km}(t,s)): t\in I_k \text{ and } s\in J_m \right\}\\
			=& \bigcup_{k=1}^{N}\bigcup_{m=1}^{M} Gr(h_{km})
			= \bigcup_{k=1}^{N}\bigcup_{m=1}^{M} W_{km}\left(A_{R_{k}Q_{m}}\right)= \boldsymbol{\Omega}(A).
		\end{align*}
		This implies that $h(t,s)=f_{[\boldsymbol{\Omega}(A)]}(t,s)$ for all $t\in I$ and $s\in J$ such that $\boldsymbol{\Omega}(A)\in \mathbf{G}_{\Delta}(I\times J, \boldsymbol{Y})$ and
		\[f_{[\boldsymbol{\Omega}(A)]}(t,s)=h_{km}(t,s)=P_{km}\left(U_{k}^{-1}(t), V_{m}^{-1}(s), f_{[A]}\left(U_{k}^{-1}(t),V_{m}^{-1}(s)\right)\right)\]
		for all $t\in I_k$ and $s\in J_m$, where $k\in \Sigma_{N}$ and $m\in \Sigma_{M}$. This completes the proof.
	\end{proof}

	\begin{theorem}\label{Th1}
		$\boldsymbol\Omega$ as defined in Equation \eqref{omega}, is a contraction map on $\left(\mathbf{G}_{\Delta}(I\times J, \boldsymbol{Y}), d_{\mathbf{G}}\right)$.
	\end{theorem}
	
	\begin{proof}
		For $A_1, A_2 \in \mathbf{G}_{\Delta}(I\times J, \boldsymbol{Y})$, we have
		\begin{align*}
			& d_{\mathbf{G}}\left(\boldsymbol\Omega(A_1), \boldsymbol\Omega(A_2)\right)\\
			=&\max_{(t,s)\in I\times J}d_{\boldsymbol{Y}}\left(f_{[\boldsymbol\Omega(A_1)]}(t,s), f_{[\boldsymbol\Omega(A_2)]}(t,s)\right)\\
			=&\max_{\substack{k=1,\ldots,N\\m=1,\ldots,M}}\max_{(t,s)\in I_k\times J_m}d_{\boldsymbol{Y}}\left(P_{km}\left(U_{k}^{-1}(t),V_{m}^{-1}(s), f_{[A_1]}\left(U_{k}^{-1}(t),V_{m}^{-1}(s)\right)\right),\right.\\ & \hspace{5cm} \left.P_{km}\left(U_{k}^{-1}(t),V_{m}^{-1}(s), f_{[A_2]}\left(U_{k}^{-1}(t),V_{m}^{-1}(s)\right)\right)\right)\\
			=&\max_{\substack{k=1,\ldots,N\\m=1,\ldots,M}}\max_{(t^*,s^*)\in R_k\times Q_m}d_{\boldsymbol{Y}}\left(P_{km}\left(t^*,s^*,f_{[A_1]}(t^*,s^*)\right), P_{km}\left(t^*,s^*,f_{[A_2]}(t^*,s^*)\right)\right)\\
			\leq&\max_{\substack{k=1,\ldots, N\\m=1,\ldots,M}}\max_{(t^*,s^*)\in R_k\times Q_m}\eta_{km}d_{\boldsymbol{Y}}\left(f_{[A_1]}(t^*,s^*),f_{[A_2]}(t^*,s^*)\right)\\
			\leq&\eta \max_{(t^*,s^*)\in \mathfrak{D}}d_{\boldsymbol{Y}}\left(f_{[A_1]}(t^*,s^*), f_{[A_2]}(t^*,s^*)\right),
		\end{align*}
		where $\eta=\max_{\substack{k=1,\ldots
				,N\\m=1,\ldots,M}}\{\eta_{km}\}$  and  $\mathfrak{D}=\bigcup_{k=1}^{N}\bigcup_{m=1}^{M}(R_k\times Q_m)$. Since $\mathfrak{D}\subseteq I\times J$ and $0<\eta<1$, then we have
		\begin{align*}
			d_{\mathbf{G}}\left(\boldsymbol\Omega(A_1), \boldsymbol\Omega(A_2)\right)&\leq \eta \max_{(t^*,s^*)\in I\times J}d_{\boldsymbol{Y}}\left(f_{[A_1]}(t^*,s^*), f_{[A_2]}(t^*,s^*)\right)=~\eta d_{\mathbf{G}}(A_1,A_2).
		\end{align*}
		This completes the proof.
	\end{proof}
	
	Now using Theorem \ref{Th1} and Banach contraction principle, $\boldsymbol\Omega$ has a fixed point $\mathbb{P}$ in $\mathbf{G}_{\Delta}(I\times J, \boldsymbol{Y})$ and $f_{[\mathbb{P}]}\in \mathcal{C}_{\Delta}(I\times J, \boldsymbol{Y})$ will be a fractal interpolation function corresponding to $\Delta$. 
	
	Define a function $T:\mathcal{C}_{\Delta}(I\times J, \boldsymbol{Y})\rightarrow \mathcal{C}_{\Delta}(I\times J, \boldsymbol{Y})$ such that 
	\begin{equation}\label{Rb1}
		T(h)=f_{[\boldsymbol\Omega(Gr(h))]}.
	\end{equation}
	This means for every $h\in \mathcal{C}_{\Delta}(I\times J, \boldsymbol{Y})$  , we have
	\begin{equation}\label{Rb2}
		T(h)(t,s)=P_{km}\left(U_{k}^{-1}(t),V_{m}^{-1}(s), h\left(U_{k}^{-1}(t),V_{m}^{-1}(s)\right)\right) \text{ for } (t,s)\in I_{k}\times J_m,
	\end{equation}
	where $k\in \Sigma_{N}$, and $m\in \Sigma_{M}$.
	
	\begin{theorem}\label{Th2}
		$T$, defined in Equation \eqref{Rb1}, is a contraction mapping on $\mathcal{C}_{\Delta}(I\times J, \boldsymbol{Y})$ and $f_{[\mathbb{P}]}$ is the fixed point of $T$ in $\mathcal{C}_{\Delta}(I\times J, \boldsymbol{Y})$. Moreover, for any $h\in \mathcal{C}_{\Delta}(I\times J, \boldsymbol{Y})$, the sequence $T^{n}(h)\to f_{[\mathbb{P}]}$ as $n\to \infty$, where $T^n(h)=(T\circ T^{n-1})(h)$ with $T^0(h)=h,~ T^{1}(h)=T(h)$.
	\end{theorem}
	
\begin{figure}[H]
		\centering
		
		\begin{subfigure}[t]{0.49\textwidth}
			\centering
			\includegraphics[width=\textwidth]{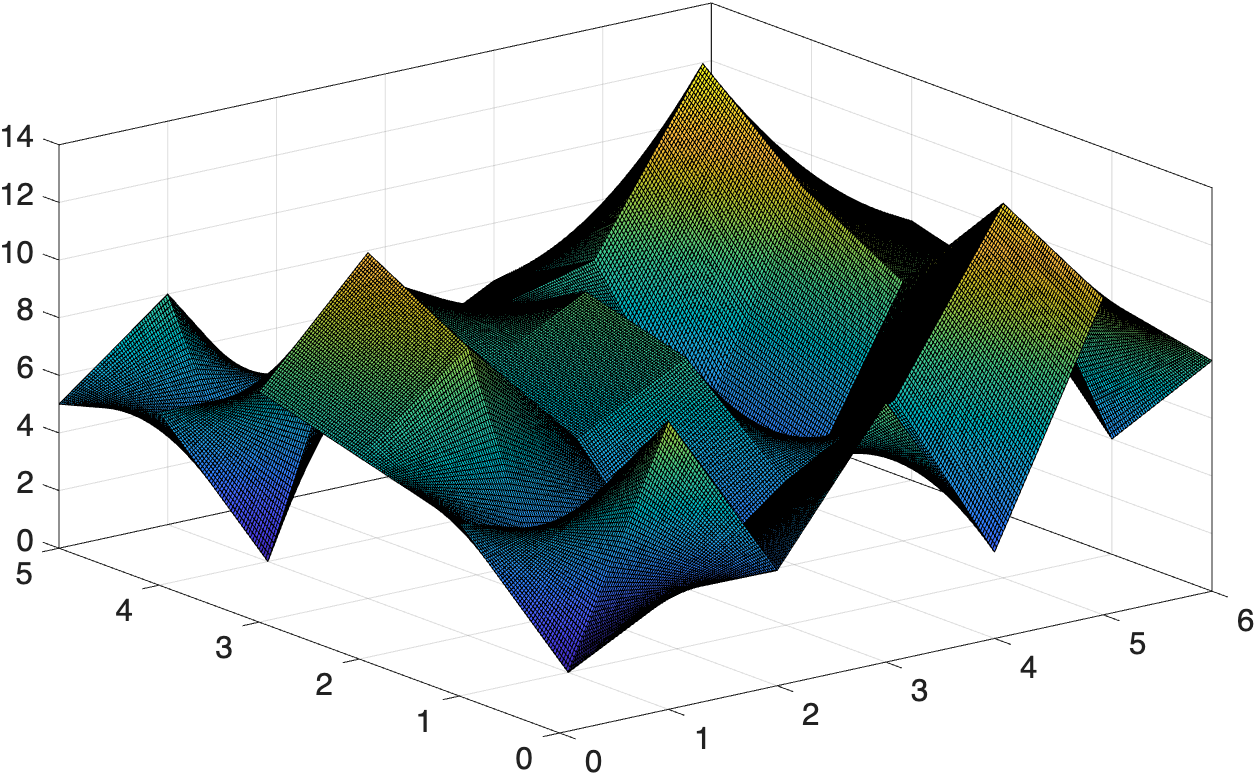}
			\caption{Bilinear interpolation surface.}
			\label{fig1}
		\end{subfigure}
		\hfill
		\begin{subfigure}[t]{0.49\textwidth}
			\centering
			\includegraphics[width=\textwidth]{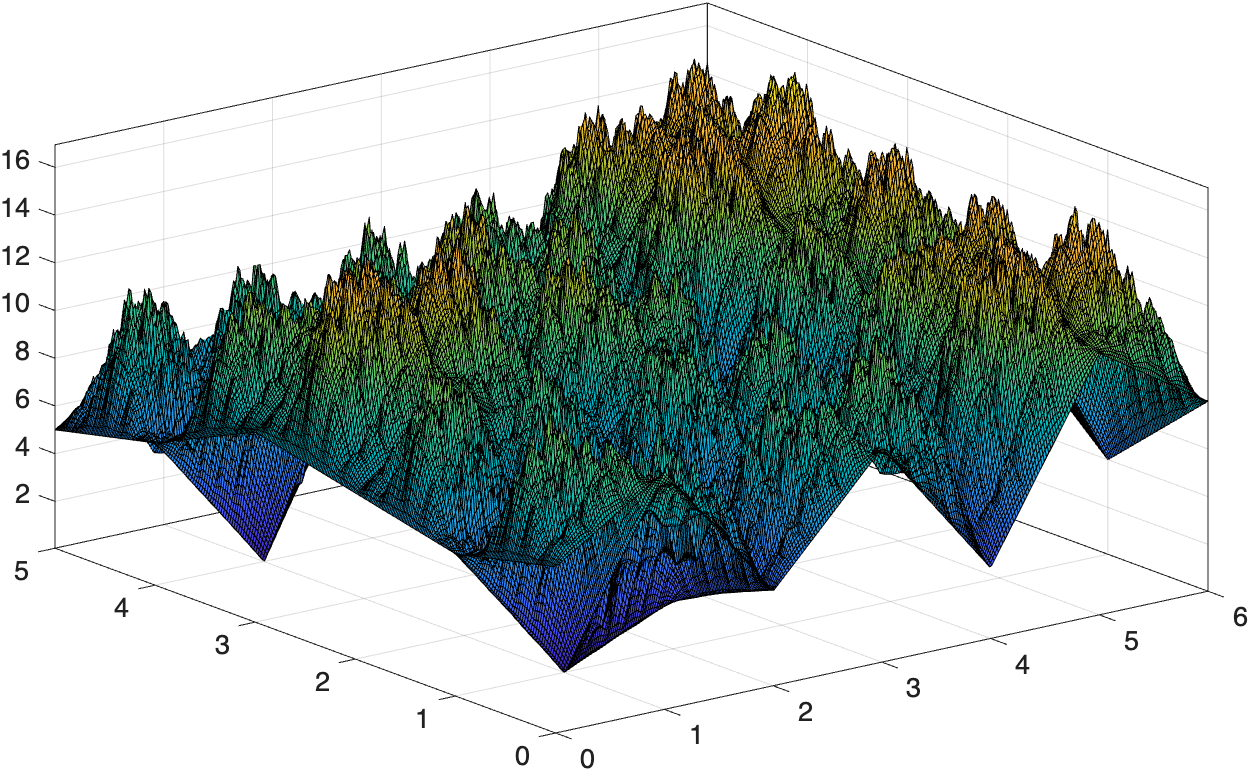}
			\caption{BFIF with $R_k=[0,6]$ and $Q_m=[0,5]$ for all
				$k=1,\ldots,6$ and $m=1,\ldots,5$.}
			\label{fig2}
		\end{subfigure}
		
	\end{figure}

\begin{figure}[H]
		\centering
		
		\begin{subfigure}[t]{0.49\textwidth}
			\centering
			\includegraphics[width=\textwidth]{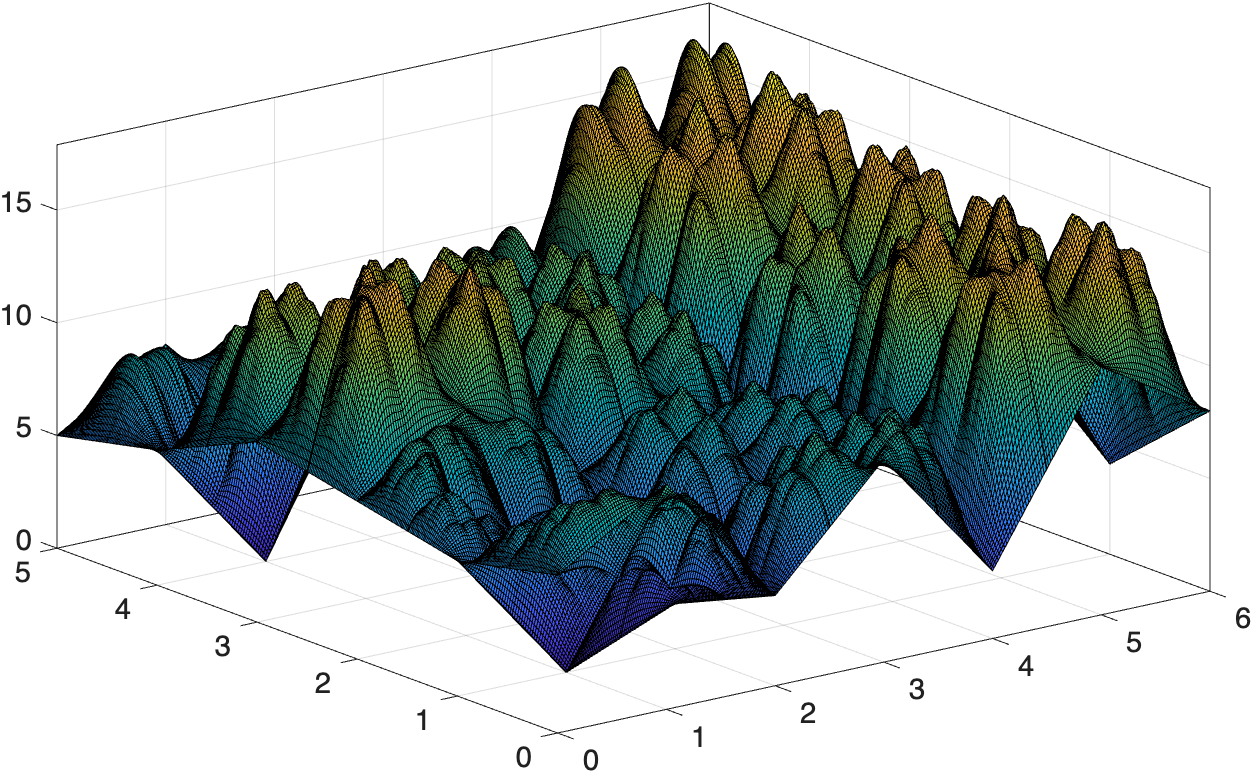}
			\caption{BFIF for $R_1 = R_2 = [0,2], \quad R_3 = R_4 = [2,4], \quad R_5 = R_6 = [4,6],$ and	$ Q_1 = Q_2 = [0,2], \quad Q_3 = Q_4 = [2,4], \quad Q_5 = [4,5].$}\label{fig3}
		\end{subfigure}
		\hfill
		\begin{subfigure}[t]{0.49\textwidth}
			\centering
			\includegraphics[width=\textwidth]{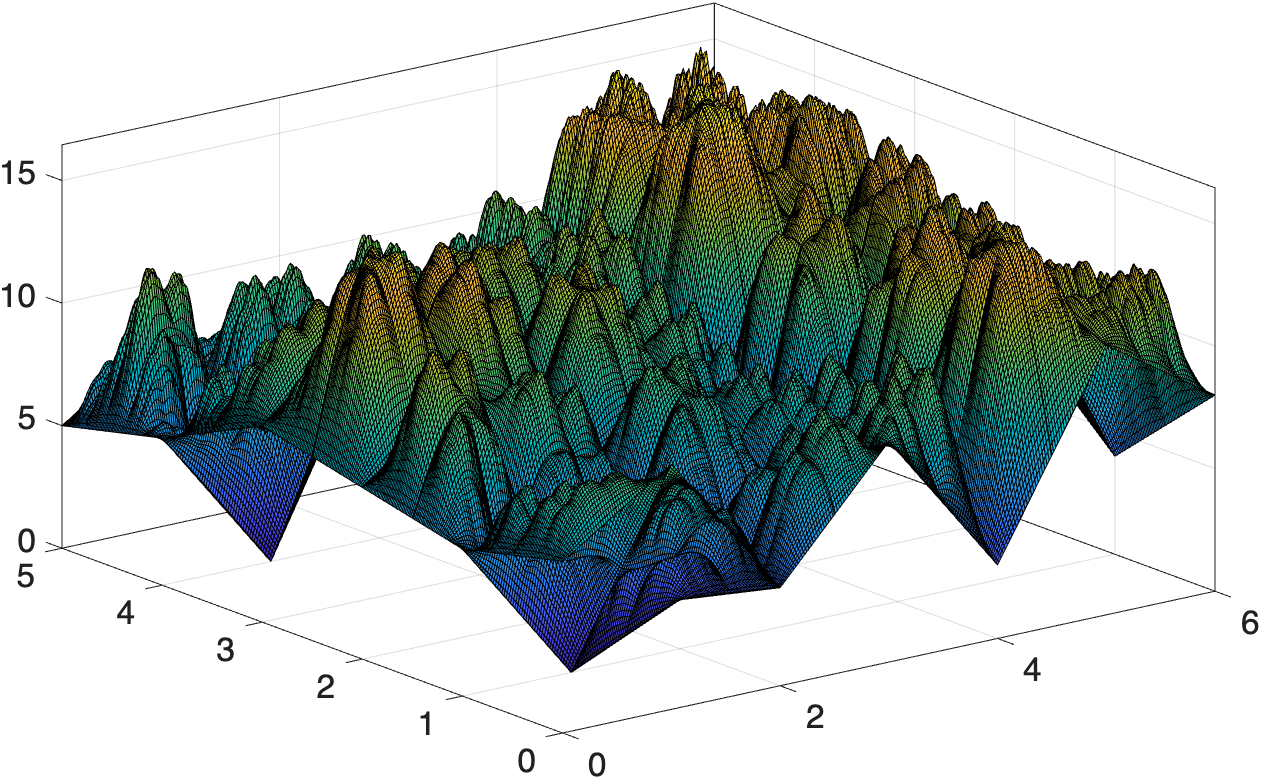}
			\caption{BFIF for $R_1 = [0,2], \quad R_2 = [1,3], \quad R_3 = [2,4], \quad R_4 = [3,5], \quad R_5 = [4,6], \quad R_6 = [0,6]$ and $Q_1 = [0,2], \quad Q_2 = [1,3], \quad Q_3 = [2,4], \quad Q_4 = [3,5], \quad Q_5 =  [0,5].$}\label{fig4}
		\end{subfigure}
	\end{figure}

\begin{figure}[H]
		\centering
		
		\begin{subfigure}[t]{0.49\textwidth}
			\centering
			\includegraphics[width=\textwidth]{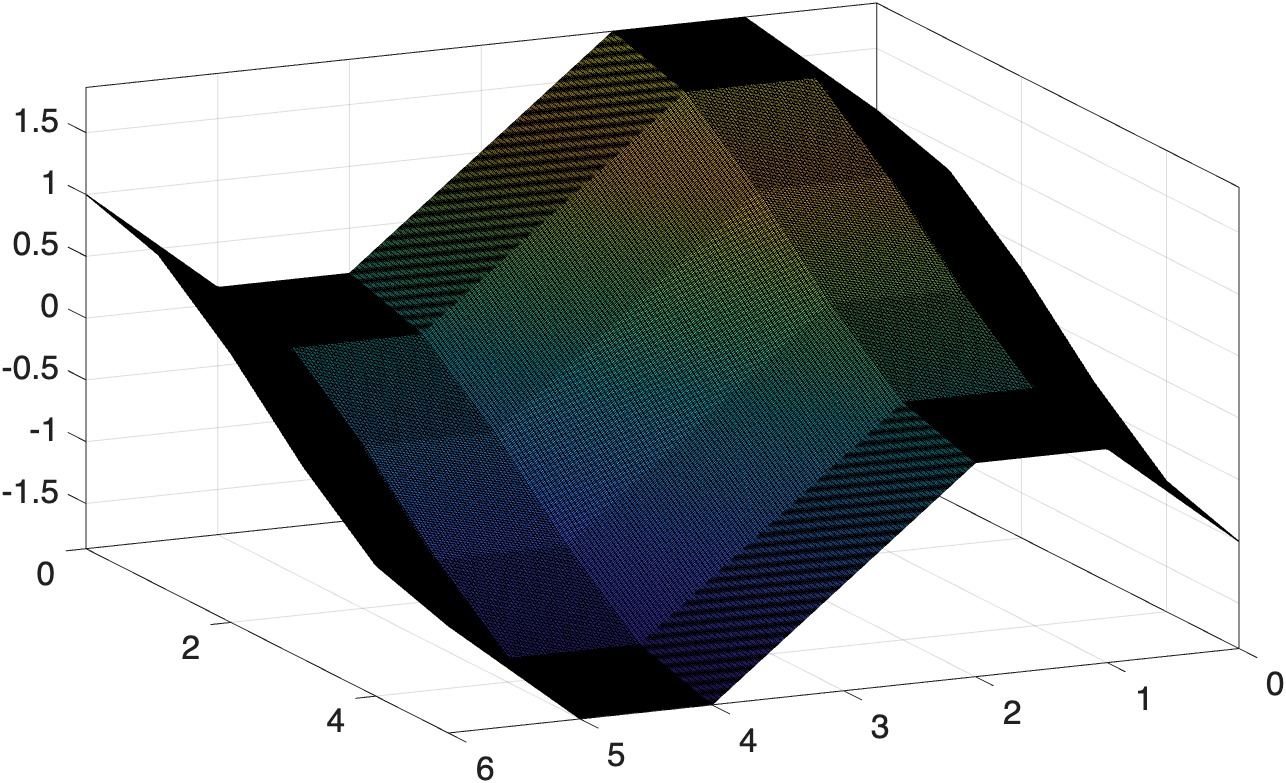}
			\caption{Bilinear interpolation surface.}\label{fig5}
		\end{subfigure}
		\hfill
		\begin{subfigure}[t]{0.49\textwidth}
			\centering
			\includegraphics[width=\textwidth]{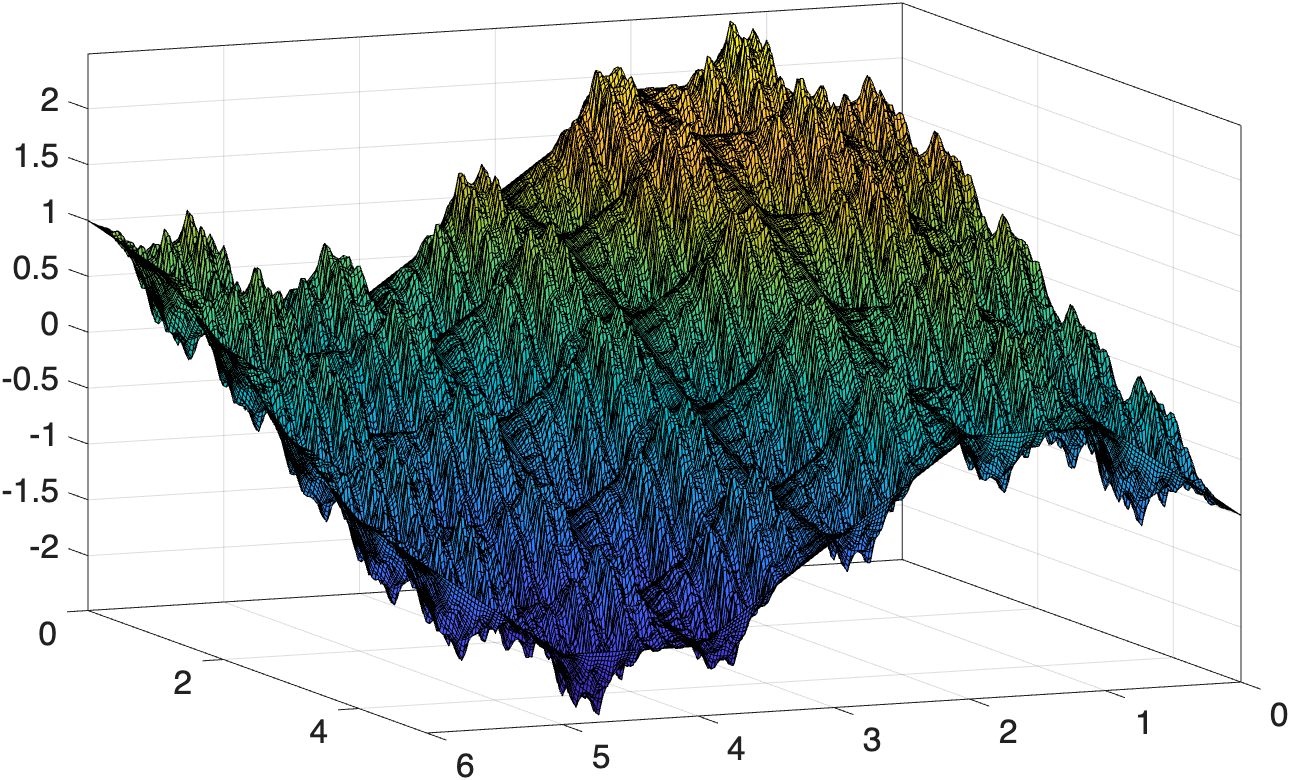}
			\caption{BFIF $R_k = [0,6]$ and $Q_m = [0,5]$ for all $k = 1,\ldots,6$ and $m = 1,\ldots,5$. }\label{fig6}
		\end{subfigure}
	\end{figure}

\begin{figure}[H]
		\centering
		
		\begin{subfigure}[t]{0.49\textwidth}
			\centering
			\includegraphics[width=\textwidth]{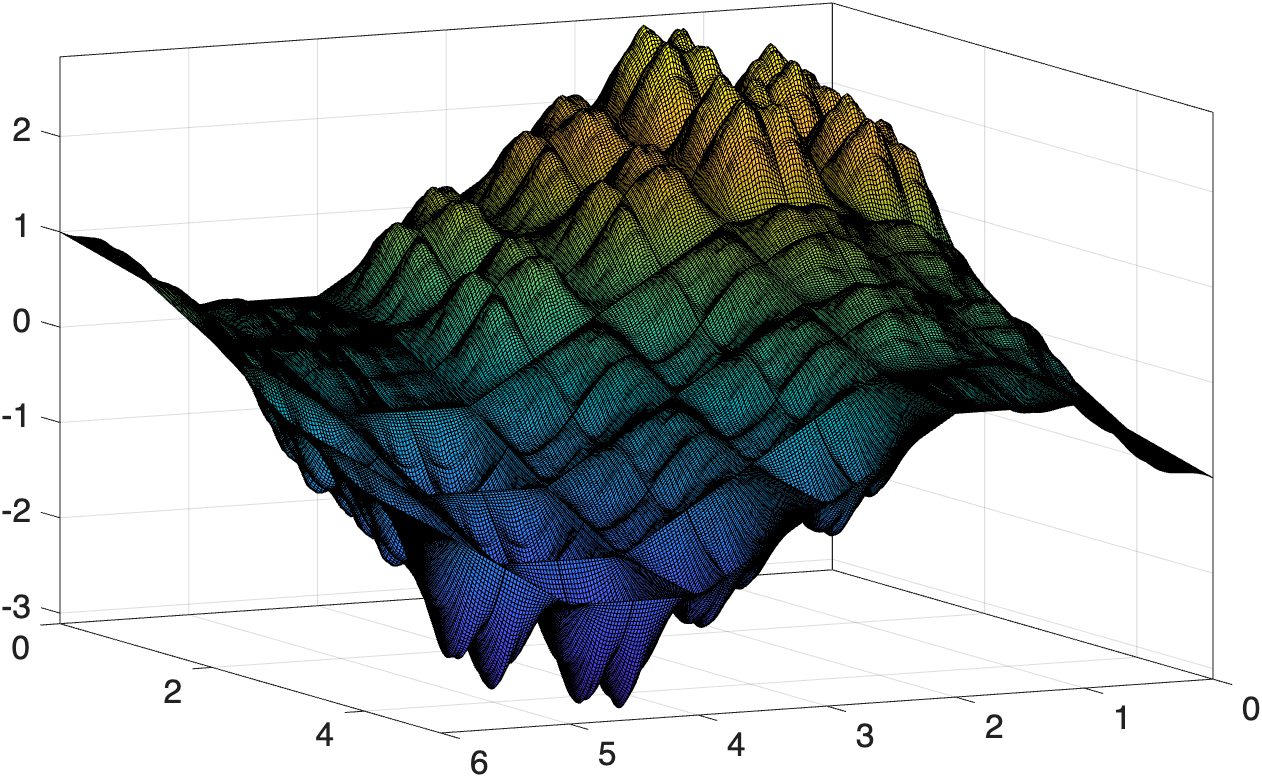}
			\caption{BFIF for $ R_1 = R_2 = [0,2], \quad R_3 = R_4 = [2,4], \quad R_5 = R_6 = [4,6],$ and $ Q_1 = Q_2 = [0,2], \quad Q_3 = Q_4 = [2,4], \quad Q_5 = [4,5].$ }\label{fig7}
		\end{subfigure}
		\hfill
		\begin{subfigure}[t]{0.49\textwidth}
			\centering
			\includegraphics[width=\textwidth]{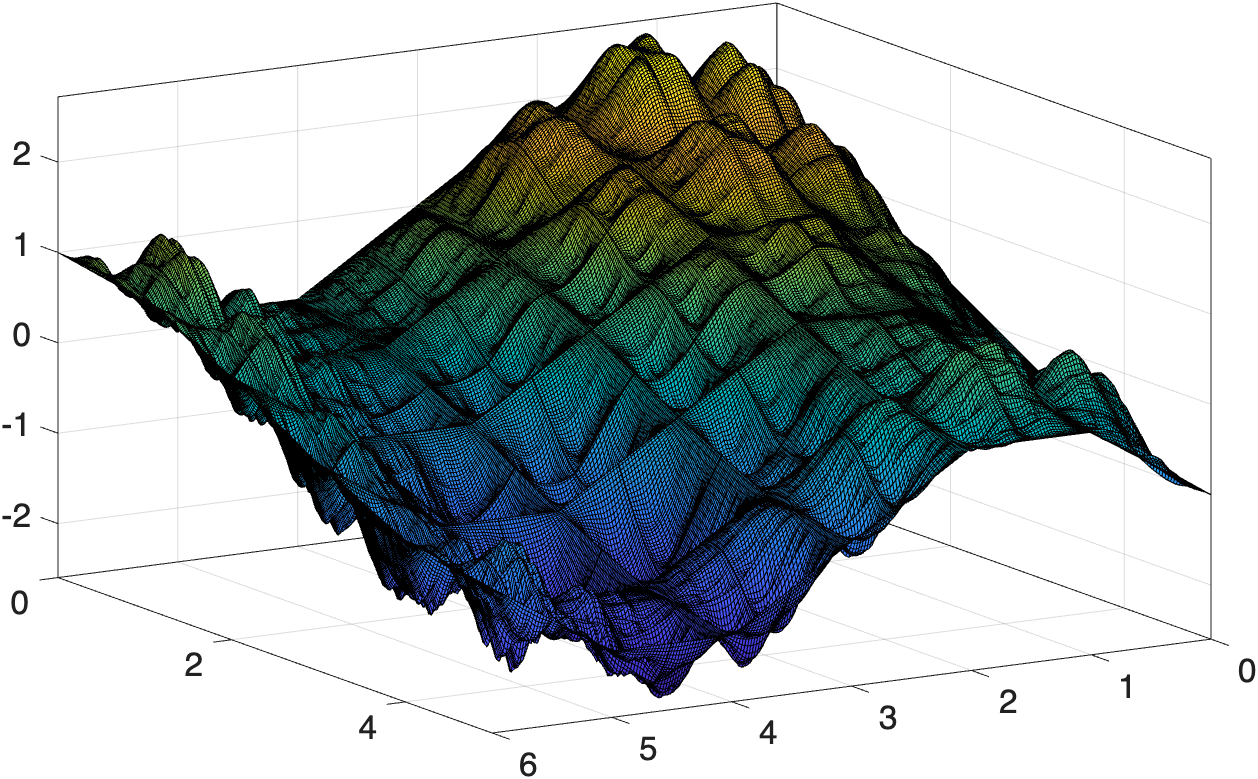}
			\caption{BFIF for $R_1 = [0,2], \quad R_2 = [1,3], \quad R_3 = [2,4], \quad R_4 = [3,5], \quad R_5 = [4,6], \quad R_6 = [0,6]$ and $Q_1 = [0,2], \quad Q_2 = [1,3], \quad Q_3 = [2,4], \quad Q_4 = [3,5], \quad Q_5 =  [0,5].$}\label{fig8}
		\end{subfigure}
	\end{figure}

	\begin{proof}
		In view of Theorem \ref{Th1}, for $h_1, h_2 \in \mathcal{C}_{\Delta}(I\times J, \boldsymbol{Y})$, we have
		\begin{align*}
			& d_{\infty}(T(h_1),T(h_2))= \max_{(t,s)\in I\times J}d_{\boldsymbol{Y}}(T(h_1)(t,s),T(h_2)(t,s))\\
			=&\max_{(t,s)\in I\times J} d_{\boldsymbol{Y}}\left(f_{[\boldsymbol\Omega(Gr(h_1))]}(t,s), f_{[\boldsymbol\Omega(Gr(h_2))]}(t,s)\right)\\
			=& \max_{\substack{k=1,\ldots,N\\m=1,\ldots,M}}\max_{(t,s)\in I_k\times J_m}d_{\boldsymbol{Y}}\left(P_{km}\left(U_{k}^{-1}(t), V_{m}^{-1}(s), h_1(U_{k}^{-1}(t), V_{m}^{-1}(s))\right),P_{km}\left(U_{k}^{-1}(t), V_{m}^{-1}(s), h_2(U_{k}^{-1}(t), V_{m}^{-1}(s))\right) \right)\\
			=& \max_{\substack{k=1,\ldots,N\\m=1,\ldots,M}}\max_{(t^*,s^*)\in R_k\times Q_m}d_{\boldsymbol{Y}}\left(P_{km}\left(t^*,s^*,h_1(t^*,s^*)\right),P_{km}\left(t^*,s^*,h_2(t^*,s^*)\right) \right)\\
			\leq& \max_{\substack{k=1,\ldots, N\\m=1,\ldots,M}}\max_{(t^*,s^*)\in R_k\times Q_m} \eta_{km} d_{\boldsymbol{Y}}\left(h_1(t^*,s^*), h_2(t^*,s^*)\right)\\
			\leq&~\eta \max_{(t^*,s^*)\in \mathfrak{D}}d_{\boldsymbol{Y}}\left(h_1(t^*,s^*), h_2(t^*,s^*)\right)\\
			\leq&~ \eta \max_{(t^*,s^*)\in I\times J} d_{\boldsymbol{Y}}\left(h_1(t^*,s^*), h_2(t^*,s^*)\right)= \eta d_{\infty}(h_1,h_2).
		\end{align*}
		Hence, $T$ is a contraction map on $\mathcal{C}_{\Delta}(I\times J, \boldsymbol{Y}).$ 
		
		Since $\mathbb{P}$ is a fixed point of $\boldsymbol{\Omega}$, therefore we have
		\begin{align*}
			&\mathbb{P}=\boldsymbol\Omega(\mathbb{P})\\
			\implies &f_{[\mathbb{P}]}(t,s) = f_{[\boldsymbol\Omega(\mathbb{P})]}(t,s)=P_{km}\left(U_{k}^{-1}(t), V_{m}^{-1}(s), f_{\mathbb{[P]}}(U_{k}^{-1}(t), V_{m}^{-1}(s))\right) \text{ for } (t,s)\in I_k\times J_m\\
			\implies & f_{[\mathbb{P}]}(t,s)= T\left(f_{[\mathbb{P}]}\right)(t,s) \text{ for } (t,s)\in I\times J\\
			\implies & f_{[\mathbb{P}]} =T\left(f_{[\mathbb{P}]}\right).
		\end{align*}
		Further, let $h$ be an element in $\mathcal{C}_{\Delta}(I\times J, \boldsymbol{Y})$, then we have
		\begin{align*}
			d_{\infty}(T^{n+1}(h), T^{n}(h)) =  d_{\infty}(T(T^{n}(h)), T(T^{n-1}(h)))
			\leq \eta d_{\infty}(T^{n}(h),T^{n-1}(h))
			\ldots
			\leq  \eta^n d_{\infty}(T(h),h).
		\end{align*}
		Now for $m,n\in \mathbb{N}$ and $m>n$, we have
		\begin{align*}
			d_{\infty}(T^{m}(h),T^{n}(h))\leq &~ \sum_{k=1}^{m-n}d_{\infty}(T^{m-(k-1)}(h),T^{m-k}(h))\\
			\leq & \sum_{k=1}^{m-n}\eta^{m-k}d_{\infty}(T(h),h)=\sum_{j=n}^{m-1} \eta^{j}d_{\infty}(T(h),h)\leq \sum_{j=n}^{\infty}\eta^{j}d_{\infty}(T(h),h).
		\end{align*}
		This implies that ${T^{n}(h)}$ is a Cauchy sequence in $\mathcal{C}_{\Delta}(I\times J, \boldsymbol{Y})$. Using completeness of $\mathcal{C}_{\Delta}(I\times J, \boldsymbol{Y})$, let $\lim_{n\to \infty}T^{n}(h)=h^*$. It is easy to notice that $f_{[\mathbb{P}]}=h^*$.
		This completes the proof.
	\end{proof}
	We present the following examples to compare the fractal interpolation surfaces obtained by using our approach than those constructed in earlier works. To highlight the differences clearly, we consider the simplest setting of linear partial fractal interpolation functions. Let the mappings $U_k$, $V_m$, and $P_{km}$ be defined by
	$U_k(t)=a_k t+b_k,\quad V_m(s)=c_m s+d_m,$	and $P_{km}(t,s,y)=\alpha_{km}y+q_{km}(t,s)$ respectively. where $q(t,s) =A_{km}t+B_{km}s + C_{km}ts + D_{km}$. Suppose that $U_k$, $V_m$, and $P_{km}$ satisfy conditions (I), (II), and (III) respectively. The exact expressions for the coefficients $a_k$, $b_k$, $c_m$, $d_m$, $A_{km}$, $B_{km}$, $C_{km}$ and $D_{km}$ can be obtained directly from the interpolation conditions. For these examples, all vertical scaling factors are chosen as $\alpha_{km}=0.5$. The figures presented below illustrate how different choices of the domains $R_k$ and $Q_m$ affect the geometric structure and partial self-similar behavior of the resulting fractal interpolation surfaces.
	\begin{example}\label{ex1}
		Let $\Delta = \{(x_i, y_j, z_{ij}) : i = 0,1,\ldots,6 \text{ and } j = 0,1,\ldots,5\}$, where  
		$x = \{0,1,2,3,4,5,6\}$ and $y = \{0,1,2,3,4,5\}$. The corresponding data values $z_{ij}$ are given by  
		\[Z =\begin{bmatrix}
			7 & 10 & 4 & 9 & 3 & 11 & 8 \\
			6 & 0  & 2 & 5 & 7 & 13 & 4 \\
			7 & 10 & 6 & 8 & 2 & 9  & 8 \\
			8 & 12 & 7 & 9 & 3 & 11 & 9 \\
			6 & 0  & 9 & 7 & 8 & 14 & 8 \\
			5 & 8  & 3 & 4 & 6 & 7  & 10
		\end{bmatrix}.\]
		The graph of the piecewise bilinear interpolation function corresponding to $\Delta$ is shown in Fig.~\ref{fig1}. The graph of the bivariate fractal interpolation surface obtained by choosing $R_k = I$ and $Q_m = J$ for all $k = 1,\ldots,6$ and $m = 1,\ldots,5$ is shown in Fig.~\ref{fig2}. The fractal interpolation surface shown in Fig.~\ref{fig3}, is obtained by choosing block domains. In particular, we take  $R_1 = R_2 = [0,2], \quad R_3 = R_4 = [2,4], \quad R_5 = R_6 = [4,6],$ and	$ Q_1 = Q_2 = [0,2], \quad Q_3 = Q_4 = [2,4], \quad Q_5 = [4,5].$
		The fractal interpolation surface, where its graph is shown in Fig.~\ref{fig4}, is constructed by choosing overlapping domains. In this case, the domains are $R_1 = [0,2], \quad R_2 = [1,3], \quad R_3 = [2,4], \quad R_4 = [3,5], \quad R_5 = [4,6], \quad R_6 = [0,6]$ and $Q_1 = [0,2], \quad Q_2 = [1,3], \quad Q_3 = [2,4], \quad Q_4 = [3,5], \quad Q_5 =  [0,5]. $
	\end{example}
	\begin{example}\label{ex2}
		Let $ \Delta^2=\{(x_i,y_j,z_{ij}) : i=0,1,\ldots,6,\; j=0,1,\ldots,5\}$, where 
		$ x=\{0,1,2,3,4,5,6\}, ~ y=\{0,1,2,3,4,5\}.$
		The data values are generated by the function
		\[z_{ij}=\sin\!\left(\frac{\pi x_i}{3}\right)+\cos\!\left(\frac{\pi y_j}{5}\right),\]
		and are given approximately by
		\[ Z^{(2)}=
		\begin{bmatrix}
			1.00 & 1.87 & 1.87 & 1.00 & 0.13 & 0.13 & 1.00 \\
			0.81 & 1.68 & 1.68 & 0.81 & -0.06 & -0.06 & 0.81 \\
			0.31 & 1.18 & 1.18 & 0.31 & -0.56 & -0.56 & 0.31 \\
			-0.31 & 0.56 & 0.56 & -0.31 & -1.18 & -1.18 & -0.31 \\
			-0.81 & 0.06 & 0.06 & -0.81 & -1.68 & -1.68 & -0.81 \\
			-1.00 & -0.13 & -0.13 & -1.00 & -1.87 & -1.87 & -1.00
		\end{bmatrix}.	\]
		The graph of the piecewise bilinear interpolation function corresponding to $\Delta^{2}$ is shown in Fig.~\ref{fig5}. The graph of the bivariate fractal interpolation surface obtained by choosing $R_k = I = [0,6], \quad k=1,\ldots,6,$ and $ Q_m = J = [0,5], \quad m=1,\ldots,5,$ is shown in Fig.~\ref{fig6}. The fractal interpolation surface, whose graph is shown in Fig.~\ref{fig7}, is constructed by choosing block domains. In particular, we take $ R_1 = R_2 = [0,2], \quad R_3 = R_4 = [2,4], \quad R_5 = R_6 = [4,6],$ and $ Q_1 = Q_2 = [0,2], \quad Q_3 = Q_4 = [2,4], \quad Q_5 = [4,5].$ The fractal interpolation surface, whose graph is shown in Fig.~\ref{fig8}, is constructed by choosing overlapping domains. In this case, the domains are $
		R_1 = [0,2], \quad R_2 = [1,3], \quad R_3 = [2,4], \quad R_4 = [3,5], \quad R_5 = [4,6], \quad R_6 = [0,6], $ and $ Q_1 = [0,2], \quad Q_2 = [1,3], \quad Q_3 = [2,4], \quad Q_4 = [3,5], \quad Q_5 = [0,5].$
	\end{example}

	\section{Partial $\alpha$-fractal function}\label{sec4}
	Throughout this section, suppose that $(\boldsymbol{Y}, \| \cdot \|_{\boldsymbol{Y}})$ is a Banach space. In this section, we discuss the univariate partial $\alpha$ fractal function. For each $k=1,\cdots, N$, let $I_k=[t_{k-1},t_k] , J_k=[t_{j(k)}, t_{l(k)}]$ such that $j(k)<l(k)$, $\Delta= \{(t_i,y_i)\in \mathbb{R}\times \mathbf{Y}: i=0,1,\cdots N\}$, where $N>1$ and $0=t_0<t_1\cdots<t_N=1$. $C[I ,\mathbf{Y}]$   denotes the set of all continuous functions defined on $I$, and $C_{\Delta}[I ,\mathbf{Y}]$ contain functions in $ C[I ,\mathbf{Y}]$ that interpolate all points in $\Delta$. Let $\alpha_k(t) $ be the variable scaling function such that $\alpha_k:J_k \to \mathbb{R}$ such that $|\alpha_k(t)|\leq s<1$ for all $k=1,\cdots, N$ and $t\in J_k$. Let $u\in C_\Delta[I,\mathbf{Y}]$ such that $u(t_i)=y_i$. Define $\Delta^*=\left\{
	(t_{j(k)},y_{j(k)}),	(t_{l(k)},y_{l(k)}):	k=1,\ldots,N	\right\}.$ Let $	b\in C_{\Delta^*}[I,\mathbf{Y}]$	be the base function, so that for every $k=1,\ldots,N$ we have $b(t_{j(k)})=y_{j(k)}$ and $b(t_{l(k)})=y_{l(k)}$. \\
	Let us define the map $ L_k:J_k\to I_k$ which is homeomorphism such that $L_k(t_{j(k)})=t_{k-1}$ and $L_k(t_{l(k)})=t_{k}$ and $M_k:J_k \times \mathbf{Y} \to \mathbf{Y}$ as \[M_k(t,y)=\alpha_k(t)y+u(L_k(t))-\alpha_k(t)b(t).\] 
	Now, define $W_k:J_k\times \mathbf{Y} \to I_k\times\mathbf{Y}$ by \[W_k(t,y)=\big(L_k(t),M_k(t,y)\big), t\in J_k, y\in\mathbf{Y}.\]
	Define an operator $T:C_\Delta[I,\mathbf{Y}]\to C_\Delta[I,\mathbf{Y}]$ by 
	\[T(f)(t)=M_k\big(L^{-1}_k(t),f(L^{-1}_k(t))\big), t\in I_k, k=1,\cdots N, f \in C_\Delta[I,\mathbf{Y}].\]
	According to \cite[Theorem 2.3]{luor2018fractal}. $T$ is a contraction mapping on $C_\Delta[I,\mathbf{Y}]$ and $f^{\alpha}$ is the fixed point of $T$, which is the graph of a continuous function denoted by $f^{\alpha} : I \to Y$, satisfying
	\[f^{\alpha}(t_k) = y_k, \text{for each} ~~k.\]
	Furthermore, $f^{\alpha}$ satisfies the self-referential equation
	\[f^{\alpha}(t) = M_k\big(L_k^{-1}(t), f^{\alpha}(L_k^{-1}(t))\big),\quad \text{for all} ~t \in I_k,\; k = 1,2,\dots,N,\] i.e.,
	\[	f^{\alpha}(t)= u(t) + \alpha_k\big(L_k^{-1}(t)\big)\big(f^{\alpha} - b\big)\big(L_k^{-1}(t)\big),\quad \text{ for all } t \in I_k,\; k = 1,2,\dots,N. \]
	Now, we are ready to prove the following Theorem.
	\begin{theorem}\label{Th3}
		Let $u \in C_\Delta[I,Y]$, $ b \in C_{\Delta^*} [I,Y]$, and for each $k = 1,2,\dots,N$, let $\alpha_k : J_k \to \mathbb{R}$ be continuous functions satisfying $\|\alpha_k\|_\infty < 1.$ Then, the operator $T : C_\Delta[I,Y] \to C_\Delta[I,Y]$ defined by \[ T(f)(t) = M_k\big(L_k^{-1}(t), f(L_k^{-1}(t))\big), \quad t \in I_k,\] is a contraction mapping. Consequently, there exists a unique function $f^\alpha \in C_\Delta[I,Y]$ such that $
		T(f^\alpha) = f^\alpha.$ Moreover, the function $f^\alpha$ satisfies the self-referential equation \[ f^\alpha(t) = u(t) + \alpha_k\big(L_k^{-1}(t)\big)\big(f^\alpha - b\big)\big(L_k^{-1}(t)\big), \quad \forall t \in I_k,\; k = 1,2,\dots,N.	\]
	\end{theorem}
	\begin{proof}
		Let $f_1, f_2 \in C_\Delta[I,Y]$. For $t \in I_k$, we have
		
		\begin{align*}
			d_\infty \big(T(f_1), T(f_2)\big)&= \max_{t\in I}d_\mathbf{Y} \big( T(f_1)(t), T(f_2)(t)\big)\\
			&=  \max_{k=1,\dots,N} \max_{t \in I_k}d_\mathbf{Y}\Big(M_k\big(L_k^{-1}(t), f_1(L_k^{-1}(t))\big), M_k\big(L_k^{-1}(t), f_2(L_k^{-1}(t))\big)\Big),\\
			&\leq  \max_{k=1,\dots,N} \max_{t^* \in J_k}d_\mathbf{Y}\Big(M_k\big(t^*, f_1(t^*)\big), M_k\big(t^*, f_2(t^*)\big)\Big),\\
			&\leq \max_{k=1,\dots,N} \|\alpha_k\|_\infty \max_{t^* \in J_k} d_\mathbf{Y}\big( f_1(t^*), f_2(t^*)\big),\\
			&< \eta \max_{t^* \in J_k} d_\mathbf{Y}\big( f_1(t^*), f_2(t^*)\big) = \eta d_\infty(f_1,f_2),
		\end{align*}
		where \[\eta = \max_{1 \le k \le N} \|\alpha_k\|_\infty < 1.\]
		Hence, $T$ is a contraction mapping on $C_\Delta[I,Y]$. Therefore, by the Banach fixed point theorem, there exists a unique function $f^\alpha \in C_\Delta[I,Y]$ such that \[T(f^\alpha) = f^\alpha.\]
		Therefore, it satisfies the self-referential equation \[ f^\alpha(t) = u(t) + \alpha_k\big(L_k^{-1}(t)\big)\big(f^\alpha - b\big)\big(L_k^{-1}(t)\big), \quad \forall t \in I_k,\; k = 1,2,\dots,N.	\]
		This completes the proof.
	\end{proof}
	An immediate consequence of the preceding Theorem \ref{Th3} is the following error estimate for the partial $\alpha$-fractal function.
	\begin{theorem}\label{Th4}
		Let $f^\alpha$ be the partial $\alpha$-fractal function associated with $u$ and $b$. Then,
		\begin{align*}
			\|f^\alpha-u \|_{\infty} \leq \frac{\eta }{1-\eta} \|\big(u - b\big)\|_{\infty},
		\end{align*}
		provided $\eta < 1$.
	\end{theorem}
	
	\subsection{Bivariate $\alpha$ Partial FIF}
	In this subsection, we discuss bivariate partial $\alpha$-fractal functions that are constructed in Section \ref{sec3} and establish their error bound. For $(t,s) \in R_k \times Q_m$, let $\alpha_{km}(t,s) $ be the variable scaling function such that $\alpha_{km}:R_k \times Q_m\to \mathbb{R}$ such that $\|\alpha_{km}\|_{R_k\times Q_m}<1.$ for all $k=1,\cdots, N$ and $m=1,\cdots,M$. Let $u\in \mathcal{C}_\Delta (I\times J, \boldsymbol{Y})$ such that $u(t_i,s_j)=y_{i,j} \text{  for all  } i=0,\ldots,N,\; j=0,\ldots,M.$ Let $	\Delta^* =
	\{ (t_{j(k)},s_{w(m)},y_{j(k),w(m)}), (t_{j(k)},s_{z(m)},y_{j(k),z(m)}),$ $(t_{l(k)},s_{w(m)},y_{l(k),w(m)}),  (t_{l(k)},s_{z(m)},y_{l(k),z(m)}) ~:~ k=1,\ldots,N,\; m=1,\ldots,M \}$ and $b\in \mathcal{C}_{\Delta^*}(I\times J, \boldsymbol{Y})$ is a base function such that $b$ agrees with $u$ at the four corner points of $R_k\times Q_m$, i.e.,
	\begin{align}
		\begin{cases}	
			b(t_{j(k)},s_{w(m)}) &= u(t_{j(k)},s_{w(m)}),\\
			b(t_{j(k)},s_{z(m)}) &= u(t_{j(k)},s_{z(m)}),\\
			b(t_{l(k)},s_{w(m)}) &= u(t_{l(k)},s_{w(m)}),\\
			b(t_{l(k)},s_{z(m)}) &= u(t_{l(k)},s_{z(m)}).
		\end{cases}
	\end{align} 
	Define $P_{km} : R_k \times Q_m \times \boldsymbol{Y} \rightarrow \boldsymbol{Y} $ by \[P_{km}(t,s,y)=\alpha_{km}(t,s)y + u\big(U_k(t),V_m(s)\big)-\alpha_{km}(t,s) b(t,s).\]
	In addition, we assume that $\{P_{km}\}$ satisfies the edge join-up conditions introduced in Section~\ref{sec3}.
	Define an operator $T:\mathcal{C}_{\Delta}(I\times J, \boldsymbol{Y})\to \mathcal{C}_{\Delta}(I\times J, \boldsymbol{Y})$ by 
	\begin{equation}
		T(f)(t,s)=P_{km}\left(U_{k}^{-1}(t),V_{m}^{-1}(s), f\left(U_{k}^{-1}(t),V_{m}^{-1}(s)\right)\right) \text{ for } (t,s)\in I_{k}\times J_m,
	\end{equation}
	Therefore, according to the Theorem \ref{Th2}, $T$ is a contraction mapping on $\mathcal{C}_{\Delta}(I\times J, \boldsymbol{Y})$. Then, $f^{\alpha}_{[\mathbb{P}]}$ is the fixed point of $T$, which is the graph of a continuous function denoted by $f^{\alpha}_{[\mathbb{P}]} : I \times J \to \boldsymbol{Y}$, satisfying
	\[f^{\alpha}_{[\mathbb{P}]}(t_i,s_j) = y_{ij}, \text{ for all } i=0,\cdots ,N \text{ and } j=0,\cdots,M.\]
	Moreover, $f^{\alpha}_{[\mathbb{P}]}$ satisfies the self-referential equation, i.e., 
	\[	f^{\alpha}_{[\mathbb{P}]}(t,s)= u(t,s) + \alpha_{km}\big(U_k^{-1}(t),V_m^{-1}(s)\big)\big(f^{\alpha}_{[\mathbb{P}]} - b\big)\big(U_k^{-1}(t),V_m^{-1}(s)\big), \text{ for all } k=1,\cdots ,N \text{ and } m=1,\cdots,M.\]
	Now, we are ready to prove the following theorem.
	\begin{theorem}\label{Th5}
		Let $u \in \mathcal{C}_{\Delta}(I\times J, \boldsymbol{Y})$, $b \in \mathcal{C}_{\Delta^*}(I\times J, \boldsymbol{Y})$ and for all $k=1,\cdots N$ and $m=1,\cdots,M$, $(t,s) \in R_k \times Q_m$, let $\alpha_{km}:R_k \times Q_m\to \mathbb{R} $ be continuous functions satisfying $\|\alpha_{km}\|_\infty < 1.$ Then, the operator $T : \mathcal{C}_{\Delta}(I\times J, \boldsymbol{Y})\to \mathcal{C}_{\Delta}(I\times J, \boldsymbol{Y})$ defined by \[ T(f)(t,s)=P_{km}\left(U_{k}^{-1}(t),V_{m}^{-1}(s), f^\alpha_{[\mathbb{P}]}\left(U_{k}^{-1}(t),V_{m}^{-1}(s)\right)\right) \text{ for } (t,s)\in I_{k}\times J_m,\] is a contraction mapping. Consequently, there exists a unique function $f^\alpha_{[\mathbb{P}]} \in \mathcal{C}_{\Delta}(I\times J, \boldsymbol{Y})$ such that $
		T(f^\alpha_{[\mathbb{P}]}) = f^\alpha_{[\mathbb{P}]}.$ Moreover, the function $f^\alpha_{[\mathbb{P}]}$ satisfies the self-referential equation \[	f^{\alpha}_{[\mathbb{P}]}(t,s)= u(t,s) + \alpha_{km}\big(U_k^{-1}(t),V_m^{-1}(s)\big)\big(f^{\alpha}_{\mathbb{[P]}} - b\big)\big(U_k^{-1}(t),V_m^{-1}(s)\big), \text{ for all } k=1,\cdots, N \text{ and } m=1,\cdots,M\]
	\end{theorem}
	\begin{proof}
		The proof is analogous to that of Theorem \ref{Th3}, hence omitted.
	\end{proof}
	An immediate consequence of the preceding Theorem \ref{Th5} is the following error estimate for the bivariate partial $\alpha$-fractal function.
	\begin{theorem}\label{errorbivariate}
		Let $f^\alpha_{[\mathbb{P}]}$ be the bivariate  partial $\alpha$-fractal function associated with $u$ and $b$. Then,
		\begin{align*}
			\|f^\alpha_{[\mathbb{P}]} - u\|_{\infty} \le \frac{\eta}{1 - \eta} \|u - b\|_{\infty},
		\end{align*}
		provided $\eta = \max_{\substack{1 \le k \le N \\ 1 \le m \le M}} \|\alpha_{km}\|_{R_k \times Q_m}<1.$
	\end{theorem}

	\section{Stability and Sensitivity Analysis of Fractal Interpolation Functions}\label{sec5}
	Throughout this section, suppose that $(\boldsymbol{Y}, \| \cdot \|_{\boldsymbol{Y}})$ is a Banach space.  In this section, we shall discuss the stability and sensitivity of FIFs, which are constructed in Section \ref{sec3}. We first consider the stability of FIFs. Let $\Delta = \{ (t_i,s_j,y_{ij}) \in I \times J \times Y : i=0,1,\cdots,N \text{ and } j=0,1,\cdots,M \}$ and $\overline{\Delta} = \{ (t_i,s_j,\overline{y_{ij}}) \in I \times J \times Y : i=0,1,\cdots,N \text{ and } j=0,1,\cdots,M \}$ be two given sets of points, where  $t_0<\cdots < t_N$ and $s_0<\cdots< s_M$. Let  $\Delta^* = \{ \left(t_{j(k)}, s_{w(m)}, y_{j(k)w(m)}\right),\left(t_{j(k)}, s_{z(m)}, y_{j(k)z(m)}\right), \left(t_{l(k)}, s_{w(m)}, y_{l(k)w(m)}\right),  \\
	\left(t_{l(k)}, s_{z(m)}, y_{l(k)z(m)}\right) : k=1,\cdots,N \text{ and } m=1,\cdots,M \}$ and 
	$\overline{\Delta}^* = \{ \left(t_{j(k)}, s_{w(m)}, \bar{y}_{j(k)w(m)}\right),\\ \left(t_{j(k)}, s_{z(m)}, \bar{y}_{j(k)z(m)}\right), \left(t_{l(k)}, s_{w(m)}, \bar{y}_{l(k)w(m)}\right),
	\left(t_{l(k)}, s_{z(m)}, \bar{y}_{l(k)z(m)}\right) : k=1,\cdots,N \text{ and } m=1,\cdots,M \}$. Define $P_{km} , \bar{P}_{km} : R_k \times Q_m \times Y \rightarrow Y $ by 
	\begin{equation}\label{eq:Pkm}
		P_{km}(t,s,y)=\alpha_{km}(t,s)y + u\big(U_k(t),V_m(s)\big)-\alpha_{km}(t,s) b(t,s),
	\end{equation}
	and
	\begin{equation}\label{eq:Pkmbar}
		\bar P_{km}(t,s,y)	=\alpha_{km}(t,s)y	+\bar u\big(U_k(t),V_m(s)\big)-\alpha_{km}(t,s)\bar b(t,s),
	\end{equation}
	where $ u\in C_{\Delta}(I\times J,Y),\quad
	\bar u\in C_{\bar\Delta}(I\times J,Y), \quad b\in C_{\Delta^{*}}(I\times J,Y),\quad \bar b\in C_{\bar\Delta^{*}}(I\times J,Y)$ and $\alpha_{km}:R_k\times Q_m\to\mathbb{R}$ is continuous with $\|\alpha_{km}\|_{\infty}\le \eta_{km}<1.$ Moreover,  $P_{km}$ and $\bar P_{km}$ satisfy the four corner join-up conditions in \eqref{cond} as well as the edge join-up conditions introduced in Section~\ref{sec3}.
	
	Let $T$ and $\bar{T}$ be operators on $C_{\Delta}(I\times J,Y)$ and $C_{\overline{\Delta}}(I\times J,Y)$, respectively, where $T$ is defined by $\eqref{eq:Pkm}$ and $\eqref{Rb2}$, and $\bar{T}$ is defined by $\eqref{eq:Pkmbar}$ and $\eqref{Rb2}$ with $P_{km}$ with being replaced by $\bar P_{km}$. By Theorem $\ref{Th2}$, $T$ and $\bar{T}$  are contraction mappings. Suppose that $f_T$ and $f_{\bar{T}}$ are the fixed points of $T$ and $\bar{T}$ in $C_{\Delta}(I\times J,Y)$ and $C_{\overline{\Delta}}(I\times J,Y)$, respectively. Then,
	\begin{align*}
		\| f_T - f_{\bar{T}} \|_{\infty} &= \max_{\substack{k=1,\ldots, N\\m=1,\ldots,M}} \max_{(t,s)\in R_k\times Q_m} \| P_{km}(t,s,f_T(t,s))- \bar P_{km}(t,s,f_{\bar{T}}(t,s)) \|_{\mathbf{Y}}\\
		&= \max_{\substack{k=1,\ldots, N\\m=1,\ldots,M}} \max_{(t,s)\in I_k\times J_m}\| \alpha_{km}\big(U^{-1}_k(t),V^{-1}_m(s)\big) \big(f_T(U^{-1}_k(t),V^{-1}_m(s))- f_{\bar{T}}(U^{-1}_k(t),V^{-1}_m(s)) \big) \\& + \big(u(t,s)-\bar u(t,s) \big) -\alpha_{km}\big(U^{-1}_k(t),V^{-1}_m(s)\big) \big(b(U^{-1}_k(t),V^{-1}_m(s))- \bar b(U^{-1}_k(t),V^{-1}_m(s)) \big)\|_{\mathbf{Y}}\\
		&\leq \eta	\| f_T - f_{\bar{T}} \|_{\infty} + 	\| u - \bar{u} \|_{\infty} + \eta	\| b -\bar b \|_{\infty}.
	\end{align*}
	Here $\eta = \max_{\substack{1\le k\le N\\ 1\le m\le M}} \eta_{km}$.
	In particular, if $u$, $\bar u$, $b$, and $\bar b$ are the
	piecewise bilinear interpolation functions corresponding to the data
	sets $\Delta$, $\overline{\Delta}$, $\Delta^{*}$, and
	$\overline{\Delta}^{*}$, respectively, then
	\[\lVert u - \bar u \rVert_{\infty}\le\max_{\substack{0\le i\le N\\ 0\le j\le M}} \lVert y_{i,j} - \bar y_{i,j} \rVert_{Y},
	\]
	and
	\[\lVert b - \bar b \rVert_{\infty} \le \max_{\substack{0\le i\le N\\ 0\le j\le M}}
	\lVert y_{i,j} - \bar y_{i,j} \rVert_{Y},
	\]
	Thus, we have the final result of stability
	\[ \| f_T - f_{\bar{T}} \|_{\infty}\leq \frac{(1+\eta)}{(1-\eta)} \max_{\substack{0\le i\le N\\ 0\le j\le M}}
	\lVert y_{i,j} - \bar y_{i,j} \rVert_{Y}. \]
	This result is analogous to \cite[Theorem 4.1] {wang2013fractal}.
	
	In the following, we consider the sensitivity of FIFs. Define $\rho_{km} : R_k \times Q_m \times Y \rightarrow Y $ by 
	\begin{equation}\label{sensitivity}
		\rho_{km}(t,s,y)=\big(\alpha_{km}(t,s) + \epsilon_{km} \beta_{km}(t,s)\big) y + q_{km}(t,s)  + \xi_{km} \phi_{km}(t,s),
	\end{equation}
	where, \[ q_{km}(t,s) = u\big(U_k(t),V_m(s)\big)-\big(\alpha_{km}(t,s)+ \epsilon_{km} \beta_{km}(t,s)\big) b(t,s),\]
	and $u\in C_{\Delta}(I\times J,Y)$ and
	$b\in C_{\Delta^{*}}(I\times J,Y)$ are such that
	$b$ agrees with $u$ at the four corner points of $R_k\times Q_m$, i.e.,
	\begin{align*}
		b(t_{j(k)},s_{w(m)}) &= u(t_{j(k)},s_{w(m)}),\\
		b(t_{j(k)},s_{z(m)}) &= u(t_{j(k)},s_{z(m)}),\\
		b(t_{l(k)},s_{w(m)}) &= u(t_{l(k)},s_{w(m)}),\\
		b(t_{l(k)},s_{z(m)}) &= u(t_{l(k)},s_{z(m)}).
	\end{align*}
	Assume that $0<|\varepsilon_{km}|<1$ and $0<|\xi_{km}|<1$.
	Let $\beta_{km}:R_k\times Q_m\to\mathbb{R}$ be continuous and
	$\phi_{km}:R_k\times Q_m\to Y$ be continuous such that
	\[\phi_{km}(t_{j(k)},s_{w(m)})=\phi_{km}(t_{j(k)},s_{z(m)})=\phi_{km}(t_{l(k)},s_{w(m)})=\phi_{km}(t_{l(k)},s_{z(m)})=0.\]
	Let
	\[\eta=\max_{\substack{1\le k\le N\\1\le m\le M}}\eta_{km},\qquad
	\varepsilon=\max_{\substack{1\le k\le N\\1\le m\le M}}|\varepsilon_{km}|,\qquad
	\xi=\max_{\substack{1\le k\le N\\1\le m\le M}}|\xi_{km}|.\]
	Define
	\[\beta=\max_{\substack{1\le k\le N\\1\le m\le M}}\max_{(t,s)\in R_k\times Q_m} |\beta_{km}(t,s)|,\qquad\phi=\max_{\substack{1\le k\le N\\1\le m\le M}} \max_{(t,s)\in R_k\times Q_m}\|\phi_{km}(t,s)\|_Y.\]
	It is easy to see that $\rho_{km}$ satisfies all the four corner join-up conditions as well as the edge join-up conditions introduced in Section~\ref{sec3}., i.e.,
	\begin{align*}
		\begin{cases}
			\rho_{km}\left(t_{j(k)}, s_{w(m)}, y_{j(k)w(m)}\right)= y_{(k-1)(m-1)},  \quad \rho_{km}\left(t_{j(k)}, s_{z(m)}, y_{j(k)z(m)}\right)= y_{(k-1)m},\\
			\rho_{km}\left(t_{l(k)}, s_{w(m)}, y_{l(k)w(m)}\right)= y_{k(m-1)}, \quad \hspace{0.6cm}\rho_{km}\left(t_{l(k)}, s_{z(m)}, y_{l(k)z(m)}\right)=y_{km},
		\end{cases}
	\end{align*}
	holds. Suppose that $0<\eta + \epsilon \beta <1$. Then, \eqref{contractive} hold with $P_{km}$ and $\eta_{km}$  being replaced by $\rho_{km}$ and $\eta + \epsilon \beta$ respectively. Let $T$ and $\bar{T}$	be two operators on $C_\Delta[I \times J, \mathbf{Y}]$, where $T$  is defined by $\eqref{eq:Pkm}$ and $\eqref{Rb2}$ and $\bar{T}$ is defined by $\eqref{sensitivity}$ and $\eqref{Rb2}$ with $P_{km}$	is being replaced by $\rho_{km}$, then by Theorem \ref{Th2}, $\bar{T}$ is a contraction mapping. Let $f_T$ and $f_{\bar{T}}$ are fixed points of $T$ and $\bar{T}$ in $C_\Delta[I \times J, \mathbf{Y}]$ respectively. Then,
	
	\begin{equation}\label{eq5.4}
		\begin{split}
			\| f_T - f_{\bar{T}} \|_{\infty} &\leq \max_{\substack{k=1,\ldots, N\\m=1,\ldots,M}} \max_{(t,s)\in R_k\times Q_m} \| P_{km}(t,s,f_T(t,s))-  \rho_{km}(t,s,f_{\bar{T}}(t,s)) \|_{\mathbf{Y}}\\
			&\leq \max_{\substack{k=1,\ldots, N\\m=1,\ldots,M}} \max_{(t,s)\in I_k\times J_m} |\alpha_{km}\big(U^{-1}_k(t),V^{-1}_m(s)\big)|\| \big(f_T(U^{-1}_k(t),V^{-1}_m(s))- f_{\bar{T}}(U^{-1}_k(t),V^{-1}_m(s)) \big) \|_{\mathbf{Y}} \\& + |\big(\epsilon_{km} \beta_{km} (U^{-1}_k(t),V^{-1}_m(s))\big)| \| f_{\bar{T}}(U^{-1}_k(t),V^{-1}_m(s))-b(U^{-1}_k(t),V^{-1}_m(s))\|_{\mathbf{Y}} \\& + |\xi_{km}|\|\phi_{km} (U^{-1}_k(t),V^{-1}_m(s))\|_\mathbf{Y}\\  &\leq \eta	\| f_T - f_{\bar{T}} \|_{\infty} + 	\varepsilon \beta \|  f_{\bar{T}} - b \|_{\infty} + \xi \phi.
		\end{split}
	\end{equation}
	Using the similar argument 
	\[ \|  f_{\bar{T}} - u \|_{\infty}  \leq (\eta+\varepsilon\beta)\|  f_{\bar{T}} - b \|_{\infty}+\xi\phi.  \]
	Then,
	\begin{align*}
		\|  f_{\bar{T}} - b \|_{\infty} &\leq\|  f_{\bar{T}} - u \|_{\infty}+\|  u - b \|_{\infty}\\& \leq  (\eta+\varepsilon\beta)\|  f_{\bar{T}} - b \|_{\infty}+\xi\phi +\|  u - b \|_{\infty},
	\end{align*}
	and this implies 
	\begin{equation}
		\|  f_{\bar{T}} - b \|_{\infty}\leq\frac{\xi\phi}{1-\eta-\varepsilon\beta} + \frac{\|  u - b \|_{\infty}}{1-\eta-\varepsilon\beta}. 
	\end{equation}
	Substituting this into \eqref{eq5.4}, we have
	\begin{equation}
		\| f_T - f_{\bar{T}} \|_{\infty} \leq  \frac{\xi\phi}{1-\eta-\varepsilon\beta} + \frac{\varepsilon \beta}{(1-\eta)(1-\eta-\varepsilon\beta)}\|  u - b \|_{\infty}.  
	\end{equation}
	This result is analogous to \cite[Theorem $5.1$]{wang2013fractal}.
	
	\section{Invariant Measure}\label{sec6}
	In this section, we study the invariant probability measure associated with the partial fractal interpolation. We consider the univariate partial fractal interpolation function introduced in Section~\ref{sec4}. Let $I=[t_0,t_N],$ $I_k=[t_{k-1},t_k],$ $k=1,\ldots,N,$
	and let $J_k=[t_{j(k)},t_{\ell(k)}]\subseteq I,$ $j(k)<\ell(k).$
	Let $f:I\to Y$ denote the corresponding partial fractal interpolation function and let
	\[G=G(f)=\{(t,f(t)):t\in I\},\]
	be its graph. For any interval $A\subseteq I$, denote the part of the graph of $f$ lying over $A$ by
	\[G_A:=\{(t,f(t)):t\in A\}.\]
	In particular, we write
	\[G_k:=G_{I_k}=\{(t,f(t)):t\in I_k\}, \text{ and } G_{J_k}:=\{(t,f(t)):t\in J_k\}.\]
	Thus,
	\[G=\bigcup_{k=1}^{N}G_k.\]
	
	Unlike the classical fractal interpolation \cite{barnsley1986fractal} setting, where the associated mappings act on the whole graph, and hence a single invariant probability measure $\mu$ on the graph $G$ satisfies
	\[\mu=\sum_{k=1}^{N}p_k\,\mu\circ W_k^{-1}.\]
	However, in PFIFs, the map
	\[W_k:J_k\times Y\longrightarrow I_k\times Y\]
	uses the graph over $J_k$ to construct the graph over $I_k$, and the interval $J_k$ may vary with $k$. Thus, the maps $W_k$ do not form an ordinary IFS on a common space. This leads naturally to a graph-directed formulation in which the invariant probability measure is described by a family of measures on the component graphs, coupled through the  partial self-similar relations. Define the set of source indices for the $k$th target component by
	\[\mathcal A_k=\{r\in\{1,\dots,N\}:I_r\subseteq J_k\}.\]
	Since $J_k=[t_{j(k)},t_{l(k)}],$
	we have $\mathcal A_k=\{j(k)+1,j(k)+2,\dots,l(k)\}.$
	Consequently,
	\[G_{J_k}=\{(t,f(t)):t\in J_k\}=\bigcup_{r\in\mathcal A_k}G_r.\]
	The maps $W_k:J_k\times \boldsymbol{Y}\longrightarrow I_k\times\boldsymbol{Y}$ such that 
	\[W_k(t,f(t))=\bigl(L_k(t),M_k(t,f(t))\bigr)=\bigl(L_k(t),f(L_k(t))\bigr).\]
	Since $L_k(J_k)=I_k$, it follows that
	\[G_k=W_k(G_{J_k})=\bigcup_{r\in\mathcal A_k}W_k(G_r).\]
	
	The decomposition of $G_{J_k}$ into the component graphs $G_r$, $r\in\mathcal A_k$, identifies precisely those components that are used in constructing $G_k$. Thus, the partial self-similar relation associated with $W_k$ can be described componentwise. For every $r\in\mathcal A_k$, define the restricted map $W_{kr} =W_k|_{G_r}:G_r\longrightarrow G_k.$ 
	
	Thus, the family
	\[\{W_{kr}:G_r\longrightarrow G_k:~ r\in\mathcal A_k,\ k=1,\ldots,N\},\]
	forms a graph-directed system whose vertices are the component graphs
	$G_1,\ldots,G_N$. A relation
	\[G_r\longrightarrow G_k\]
	indicates that the graph over $I_r$ is used to construct the graph over $I_k$. Consequently, the graph-directed system records the partial self-similar relations between the component graphs. 
	
	Now for every pair $(k,r)$, assume that there exists $s_{kr}\in[0,1)$ such that
	\[d_X\bigl(W_{kr}(z),W_{kr}(z')\bigr)\leq s_{kr}d_X(z,z'),	\qquad z,z'\in G_r,\]
	set $s=\max_{\substack{1\leq k\leq N\\r\in\mathcal A_k}}s_{kr}<1$.
	
	We now prove the existence and uniqueness of the corresponding invariant measure. Instead of considering a single invariant probability measure on the whole graph, we consider a vector of Borel probability measures
	\[\boldsymbol{\mu} =(\mu_1,\ldots,\mu_N),\]
	where $\mu_k$ is a probability measure supported on the component graph $G_k$. These measures are not considered independently. The measure $\mu_k$ is determined by the measures on those component graphs $G_r$ for which $r\in\mathcal A_k$, that is, by the component graphs contained in $G_{J_k}$ that are used to construct $G_k$. Hence, the invariant probability measures are coupled according to the partial self-similar relations.
	
	More precisely, we need a vector
	$\boldsymbol{\mu}^{*}=(\mu_1^{*},\ldots,\mu_N^{*})$
	satisfying
	\[\mu_k^{*}=\sum_{r\in\mathcal A_k}p_{kr}\,\mu_r^{*}\circ W_{kr}^{-1}, \qquad k=1,\ldots,N.\]
	\begin{theorem}\label{Th6}
		Let $\{W_k\}_{k=1}^N$ be the partial fractal interpolation system defined by $W_k(t,y)=\bigl(L_k(t),M_k(t,y)\bigr)$, and $W_k:J_k\times \boldsymbol{Y} \longrightarrow I_k\times \boldsymbol{Y}$.
		For each $k$ define $	\mathcal A_k=\{r\in\{1,\ldots,N\}:I_r\subseteq J_k\}$ 
		and for $r\in\mathcal A_k$, let $	W_{kr}=W_k|_{G_r}:G_r\longrightarrow G_k.$
		Assume that each $W_{kr}$ is a contraction mapping. Let $p_{kr}>0$ be the probability associated with the map $W_{kr}$ such that $	\sum_{r\in\mathcal A_k}p_{kr}=1 \text{ for }  k=1,\ldots,N.$
		Then, there exists a unique vector of Borel probability measures $ \boldsymbol{\mu}^{*} =(\mu_1^{*},\mu_2^{*},\ldots,\mu_N^{*})$
		such that
		\[\mu_k^{*}	= \sum_{r\in\mathcal A_k} p_{kr}\mu_r^{*}\circ W_{kr}^{-1},
		\qquad k=1,\ldots,N.\]
		Moreover, $\operatorname{supp}(\mu_k^{*})=G_k,~ k=1,\ldots,N,$
		and consequently, $	\bigcup_{k=1}^N\operatorname{supp}(\mu_k^{*})=G.$
	\end{theorem}
	\begin{proof}
		Let $\mathcal P(G_k)$ denote the space of all Borel probability measures on $G_k$ for $k=1,\ldots, N$, and define the Monge--Kantorovich metric by
		\[d_{MK}(\mu,\nu)=\sup_{\operatorname{Lip}(\varphi)\leq1}\left| \int_{G_k}\varphi\,d\mu	- \int_{G_k}\varphi\,d\nu\right|,\]
		where $\mu,\nu\in\mathcal P(G_k)$. Then, $ \bigl(\mathcal P(G_k),d_{MK}\bigr)$ is a complete metric space. Now consider the product space $ \mathfrak P=  \prod_{k=1}^N\mathcal P(G_k)$ such that for $\boldsymbol{\mu}=(\mu_1,\dots,\mu_N), \boldsymbol{\nu}=(\nu_1,\dots,\nu_N)$ in $\mathfrak P$, define
		\[D(\boldsymbol{\mu},\boldsymbol{\nu})=\max_{1\leq k\leq N}d_{MK}(\mu_k,\nu_k). \]
		Then, $(\mathfrak P, D)$ is a complete metric space. Define the graph-directed measure operator $ \mathcal M:\mathfrak P\longrightarrow\mathfrak P$ by
		\begin{equation}\label{measureoperator}
			(\mathcal M\boldsymbol{\mu})_k= \sum_{r\in\mathcal A_k} p_{kr}\mu_r\circ W_{kr}^{-1}, \qquad k=1,\ldots,N.
		\end{equation}
		It is easy to see that $\mathcal M$ is well defined. We next prove that $\mathcal M$ is a contraction. Let $\boldsymbol{\mu}= (\mu_1,\dots,\mu_N),  \boldsymbol{\nu} =	(\nu_1,\dots,\nu_N) \in \mathfrak P$. Fix $k\in\{1,\dots,N\}$ and let
		$\varphi:G_k\to\mathbb R$ be a Lipschitz function satisfying $ \operatorname{Lip}(\varphi)\leq1$. Using the definition of the pushforward measure, we obtain
		\[\int_{G_k}\varphi\,d(\mathcal M\boldsymbol{\mu})_k-\int_{G_k}\varphi\,d(\mathcal M\boldsymbol{\nu})_k = \sum_{r\in\mathcal A_k}p_{kr} \left( 	\int_{G_r}\varphi\circ W_{kr}\,d\mu_r - \int_{G_r}\varphi\circ W_{kr}\,d\nu_r \right).\]
		Since $\operatorname{Lip}(\varphi)\leq1$ and $W_{kr}$ is a
		contraction with constant $s_{kr}$. Then,
		\[\left|\int_{G_r}\varphi\circ W_{kr}\,d\mu_r-\int_{G_r}\varphi\circ W_{kr}\,d\nu_r\right| \leq s_{kr}d_{MK}(\mu_r,\nu_r).\]
		Therefore,
		\begin{align*}
			\left|\int_{G_k}\varphi\,d(\mathcal M\boldsymbol{\mu})_k-\int_{G_k}\varphi\,d(\mathcal M\boldsymbol{\nu})_k\right|&\quad\leq
			\sum_{r\in\mathcal A_k}p_{kr}s_{kr}d_{MK}(\mu_r,\nu_r)\\
			&\quad\leq s\sum_{r\in\mathcal A_k} p_{kr}D(\boldsymbol{\mu},\boldsymbol{\nu})\\
			&\quad=sD(\boldsymbol{\mu},\boldsymbol{\nu}).
		\end{align*}
		Taking the supremum over all $\varphi$ satisfying
		$\operatorname{Lip}(\varphi)\leq1$, we obtain
		\[d_{MK}\bigl((\mathcal M\boldsymbol{\mu})_k,(\mathcal M\boldsymbol{\nu})_k\bigr) \leq sD(\boldsymbol{\mu},\boldsymbol{\nu}).	\]
		Taking maximum over $k=1,\dots,N$ yields
		\[D\bigl(\mathcal M\boldsymbol{\mu},\mathcal M\boldsymbol{\nu}\bigr) \leq sD(\boldsymbol{\mu},\boldsymbol{\nu}). \]
		Since $s<1$, the operator $\mathcal M$ is a contraction on the complete metric space $(\mathfrak P,D)$. Therefore, by the Banach fixed point theorem, there exists a unique vector $\boldsymbol{\mu}^*= (\mu_1^*,\dots,\mu_N^*) \in\mathfrak P$ such that
		\[\mathcal M\boldsymbol{\mu}^*=\boldsymbol{\mu}^*.\]
		Now to identify the supports of the measures $\mu_k^*$. 
		Let $S_k:=\operatorname{supp}(\mu_k^*),~ k=1,\ldots,N$. Since 
		\[\mu_k^*=\sum_{r\in \mathcal A_k}p_{kr}\mu_r^{*}\circ W_{kr}^{-1},\]
		and the positivity of every \(p_{kr}\), the support of \(\mu_k^*\) is the union of the supports of the corresponding pushed-forward measures. Moreover, since each \(W_{kr}\) is continuous and \(G_r\) is compact, we have
		\[\operatorname{supp}\bigl(\mu_r^{*}\circ W_{kr}^{-1},\bigr)=		W_{kr}\bigl(\operatorname{supp}(\mu_r^*)\bigr).\]
		Consequently,
		\[S_k=\bigcup_{r\in \mathcal A_k}W_{kr}(S_r),\qquad k=1,\ldots,N.\]
		Let $\mathcal K(G_k)$ denote the space of nonempty compact subsets
		of $G_k$, equipped with the Hausdorff metric $h_k$, and define
		\[\mathfrak K=\prod_{k=1}^N\mathcal K(G_k),\qquad H(\boldsymbol{K},\boldsymbol{L})=\max_{1\leq k\leq N}h_k(K_k,L_k).\]
		Define $\mathcal W:\mathfrak K\to\mathfrak K$ by
		\[(\mathcal W\boldsymbol{K})_k=\bigcup_{r\in\mathcal A_k}W_{kr}(K_r).\]
		Since every $W_{kr}$ has contraction constant at most $s<1$,
		\[H(\mathcal W\boldsymbol{K},\mathcal W\boldsymbol{L})\leq sH(\boldsymbol{K},\boldsymbol{L}).\]
		Thus, $\mathcal W$ is a contraction and has a unique compact set as fixed
		point. Since both $(S_1,\ldots,S_N)$ and $(G_1,\ldots,G_N)$ are fixed
		points of $\mathcal W$, we obtain
		\[S_k=G_k,\qquad k=1,\ldots,N.\]
		Therefore, $\operatorname{supp}(\mu_k^*)=G_k,~ k=1,\ldots,N$ and hence
		\[\bigcup_{k=1}^{N}\operatorname{supp}(\mu_k^*)=\bigcup_{k=1}^{N}G_k=G.\]
		Thus, the union of the supports of the invariant component measures coincides with the entire graph \(G\). This completes the proof.
	\end{proof}
	
	\begin{remark}
		The vector form of the invariant probability measure in Theorem~\ref{Th6} is a consequence of the partial self-similar structure. In the classical fractal interpolation setting, every component graph is generated from the whole graph, and hence a single invariant probability measure $\mu$
		satisfies
		\[\mu=\sum_{k=1}^{N}p_k\,\mu\circ W_k^{-1}.\]
		In contrast, in the present setting, the graph over $I_k$ is generated from the graph over $J_k$, which is decomposed into the component graphs $G_r$, $r\in\mathcal A_k$. Therefore, the invariant probability measure associated with $G_k$ is coupled with the measures associated with these component graphs through
		\[\mu_k^{*}=\sum_{r\in\mathcal A_k}p_{kr}\,\mu_r^{*}\circ W_{kr}^{-1}.\]
		Thus, the vector $\boldsymbol{\mu}^{*} =(\mu_1^{*},\ldots,\mu_N^{*})$ describes the invariant probability measure on the component graphs while preserving the partial self-similar relations between them.
	\end{remark}
	
	\begin{remark}
		If $J_k=I$ for every $k$ and each map $W_k:G\longrightarrow G$ is a contraction, then the maps $W_1,\ldots,W_N$ form an ordinary IFS on the common space $G$. In that case, given probabilities $q_k>0$ satisfying $	\sum_{k=1}^Nq_k=1,$ the usual invariant-measure operator is $\mathcal T\mu = \sum_{k=1}^Nq_k~\mu\circ W_k^{-1}.$ Consequently, there exists a unique Borel probability measure $\mu^*$
		on $G$ satisfying $	\mu^*	=	\sum_{k=1}^N q_k ~\mu^*\circ W_k^{-1} .$ Thus, the graph-directed construction reduces to the usual invariant-measure
		setting for an ordinary IFS.
	\end{remark}
	
	\section{Box Dimension for Univariate PFIF}\label{sec7}
	In this section, we investigate the box dimension of univariate PFIFs under certain additional conditions. Under these conditions, the resulting PFIF framework coincides with the self-affine RFIF setting studied by Ruan et al. \cite{ruan2021existence}. However, the proof presented here is developed directly in the PFIF setting and follows a different approach than that used in \cite{ruan2021existence}. Let $ \Delta=\{(t_i,y_i):i=0,1,\dots,N\}$ such that $t_0<t_1<\cdots<t_N$ and $I=[t_0,t_N], I_k=[t_{k-1},t_k]$ for $k=1,\dots,N$ with $J_k=[t_{j(k)},t_{l(k)}]$ such that $j(k)<l(k)$ and maps $ W_k(t,y)=(L_k(t),M_k(t,y)),$ where $L_k:J_k\to I_k$
	is affine such that	
	\begin{equation*}\label{eq:a-def}
		0<a_k:=\frac{|I_k|}{|J_k|}<1.
	\end{equation*}
	Assume also that
	\begin{equation}\label{eq:M-affine}
		M_k(t,y)=c_k t+\alpha_k y+e_k,
		\qquad 0<|\alpha_k|<1.
	\end{equation}
	$M_k:J_k\times Y\to Y$ is continuous and a contraction with respect to $y$. The resulting operator on interpolation graphs has a unique fixed graph $G(f)$, and the fixed function satisfies
	\begin{equation}\label{eq:luor-functional} 
		f(t)=M_k\bigl(L_k^{-1}(t),f(L_k^{-1}(t))\bigr),\qquad t\in I_k.
	\end{equation}
	Equivalently,
	\begin{equation}\label{eq:local-graph}
		G_k=G(f)|_{I_k}=W_k\bigl(G(f)|_{J_k}\bigr).
	\end{equation}
	Throughout this section, we take $Y=\mathbb{R}$.


	\begin{condition}\label{markov}
		For each $k=1,\ldots,N$, since $J_k$ is a union of some of the intervals $I_1,\ldots,I_N$. Define the connection matrix
		\[C=(C_{km})_{N\times N}\]
		by
		\[C_{km}=
		\begin{cases}
			1, & I_m\subseteq J_k,\\
			0, & I_m\not\subseteq J_k.
		\end{cases}\]
		Denote
		\[J_k=\bigcup_{m:C_{km}=1} I_m.\]
		Assume that the directed graph determined by $C$ is irreducible.
	\end{condition}
	
	From \eqref{eq:local-graph} and Condition \ref{markov}, we get
	\begin{equation}\label{eq:graph-directed}
		G_k =\bigcup_{m:C_{km}=1} W_k(G_m),
		\quad k=1,\dots,N.
	\end{equation}
	Thus, partial self-similarity becomes a graph-directed self-affine structure.
	\begin{lemma}[Perron-Frobenius theorem]\label{perron}
		Let $A = (a_{ij})_{n\times n}$ be an irreducible nonnegative matrix
		\begin{itemize}
			\item  The spectral radius of $A$, denoted by $\rho(A)$, is an eigenvalue of $A$ and has a positive eigenvector. 
			\item $\rho(A)$ strictly increases as any $a_{ij}$ strictly increases.
		\end{itemize}
		
	\end{lemma}

	\begin{definition}
		For $s\in[1,2]$, define the matrix
		\begin{equation*}
			B(s)=\bigl(B_{km}(s)\bigr)_{k,m=1}^{N},
			\qquad
			B_{km}(s)=C_{km}|\alpha_k|a_k^{s-1}.
		\end{equation*}
		Let $\rho(B(s))$ denote its spectral radius. Since $0<a_k<1$, the function $s\mapsto \rho(B(s))$ is continuous and strictly decreasing. Define
		\begin{equation}\label{eq:Dstar}
			D_*=
			\begin{cases}
				1,& \rho(B(1))\leq 1,\\[4pt]
				\text{the unique }d\in(1,2]\text{ satisfying }\rho(B(d))=1,
				& \rho(B(1))>1.
			\end{cases}
		\end{equation}
	\end{definition}
	
	\begin{remark} 
		The existence of a unique $d$ in the definition of $D_*$ can be established using the intermediate value theorem, with $d\in (1,2]$. Indeed, let $\ell=(|I_1|,\ldots,|I_N|)^T$. For each \(k\), 
		\[ (B(2)\ell)_k = \sum_{m=1}^N B_{km}(2)|I_m| = |\alpha_k|a_k\sum_{m=1}^N C_{km}|I_m|. \] 
		Since \(C_{km}=1\), when \(I_m\subseteq J_k\), we have $\sum_{m=1}^N C_{km}|I_m|=|J_k|.$ Therefore, 
		\[ (B(2)\ell)_k = |\alpha_k|a_k|J_k| = |\alpha_k||I_k| = |\alpha_k|\ell_k < \ell_k. \] Hence, \(B(2)\ell<\ell\). By the Perron--Frobenius theorem applied to \(B(2)^T\), there exists \(w>0\) such that 
		\[ w^TB(2)=\rho(B(2))w^T. \] 
		Consequently,
		\[\rho(B(2))\,w^T\ell = w^TB(2)\ell < w^T\ell.\] 
		Since \(w^T\ell>0\), it follows that $\rho(B(2))<1.$
	\end{remark}
	
	For an interval $U\subset I$, define $ \osc(f,U)=\sup\{|f(x)-f(y)|:x,y\in U\}$. For $\delta>0$ and $\tau\in[0,\delta)$, define
	\[\mathcal P_{\delta,\tau}(U)=\left\{U\cap[\tau+n\delta,\tau+(n+1)\delta):n\in\mathbb Z,\;U\cap[\tau+n\delta,\tau+(n+1)\delta)\neq\varnothing\right\}.
	\]
	Thus, $\mathcal P_{\delta,\tau}(U)$ is the partition of $U$ induced by the translated $\delta$-mesh. Define the total oscillation of $f$ over $U$ at scale $\delta$ by
	\[O_U(\delta)=\sup_{0<\tau<\delta}\sum_{Q\in\mathcal P_{\delta,\tau}(U)}\osc(f,Q).\]
	Therefore, for interval $I_k$
	\begin{equation}\label{eq:Ok}
		O_k(\delta)=\sup_{0<\tau<\delta}\sum_{Q\in\mathcal P_{\delta,\tau}(I_k)}\osc(f,Q).
	\end{equation}
	The following results will be used later.
	
	\begin{lemma}\label{lem:osc-cover}
		There are constants $A_1,A_2>0$, independent of small $\delta$, such that
		\begin{equation*}
			N_\delta(G_k)
			\leq A_1\left(\frac{|I_k|}{\delta}+\frac{O_k(\delta)}{\delta}+1\right),
		\end{equation*}
		and
		\begin{equation*}
			N_\delta(G_k)
			\geq A_2\frac{O_k(\delta)}{\delta},
		\end{equation*}
		whenever the oscillations are measured over a disjoint mesh.
	\end{lemma}
	\begin{proof}
		Over an interval $Q$ of length $\delta$, the graph above $Q$ lies inside $Q \times [\min_Q f, \max_Q f]$. Therefore, the number of squares of side length $\delta$ are needed to cover $G_K$. 
		\begin{align*}
			N_\delta(G_K) \leq C \big(1+\frac{\osc(f,Q)}{\delta}\big).
		\end{align*}
		Let $\mathcal P_{\delta,\tau}(I_k) =\{Q_1, Q_2, \cdots, Q_M\}$ partition $I_k$ into $M$ intervals of length $\delta$. Then, $M \leq 1+ \frac{|I_k|}{\delta}$. Therefore, $G_{Q_i}=G_k \cap (Q_i \times \mathbb{R})$ be such that $G_k =\cup_{i=1} ^M G_{Q_i}$. Therefore, 
		\begin{align*}
			N_\delta(G_k) &\leq \sum_{i=1}^{M} N_\delta(G_{Q_i}) \leq  \sum_{i=1}^{M} C \big(1+\frac{\osc(f,Q_i)}{\delta}\big)\\
			&\leq C M +C \frac{ O_k(\delta)}{\delta}\\
			&\leq A_1\left(\frac{|I_k|}{\delta}+\frac{O_k(\delta)}{\delta}+1\right).
		\end{align*}
		Now, over a disjoint interval $Q$, at least $c \frac{\osc(f,Q)}{\delta}$ squares are required to cover the graph above $Q$.i.e.,
		\begin{align*}
			N_\delta(G_k) &\geq \frac{c}{\delta} \sum_{i=1}^{M} \osc(f,Q_i) \geq A_2 \frac{O_k(\delta)}{\delta}.
		\end{align*}
		This completes the proof of the lemma.
	\end{proof}
	
	Since $L_k:J_k \to I_k$ for $t\in I_k$, write $u=L_k^{-1}(t)$. From \eqref{eq:M-affine} and \eqref{eq:luor-functional}, $ f(t)=c_kL_k^{-1}(t)+\alpha_k f(L_k^{-1}(t))+e_k.$
	Let $Q\subset I_k$ be an interval of length $\delta$. Then, $L_k^{-1}(Q)$ is an interval of length $\delta/a_k$ contained in $J_k$. Hence,
	\begin{equation*}
		\osc(f,Q)
		\leq |\alpha_k|\,\osc\bigl(f,L_k^{-1}(Q)\bigr)+C_k\frac{\delta}{a_k},
	\end{equation*}
	where $C_k=|c_k|$ up to an absolute constant. Taking the supremum over all translated $\delta$-meshes of $I_k$ gives
	\begin{equation}\label{osc-rec-upper}
		O_k(\delta)
		\leq |\alpha_k|\sum_{m:C_{km}=1}O_m(\delta/a_k)+C,
	\end{equation}
	where $C$ is independent of $\delta$. Now we establish the upper box dimension theorem.
	
	
	
	\begin{theorem}\label{thm:upper-box}
		Let f be an affine PFIF with $0  < a_k <1 $ and Condition \ref{markov} holds, then $\overdimB G(f)\leq D_*$. In particular, if $\rho(B(1))\leq1$, then $\overdimB G(f)=1.$
	\end{theorem}
	
	\begin{proof}
		The lower bound $\overdimB G(f)\geq1$ is obvious because $G(f)$ is the graph of a continuous function over an interval. Hence, it remains to prove the upper bound.
		Let $s>D_*$. Then, by the definition of $D_*$, $\rho(B(s))<1.$
		
		By Perron-Frobenius theory, since $B(s)$ is a nonnegative irreducible matrix with spectral radius less than $1$, there exists a positive vector $v=(v_1,\dots,v_N)$ and a number $\theta\in(0,1)$ such that $B(s)v = \theta v, \text{~where~} \theta= \rho(B(s))<1$. For every $k=1,\dots,N$, define
		\begin{equation}\label{Hdelta}
			H_k(\delta)=\frac{\delta^{s-1}O_k(\delta)}{v_k}.
		\end{equation}
		Such that
		\[O_m(\delta/a_k)=v_m\left(\frac{\delta}{a_k}\right)^{1-s} H_m(\delta/a_k).\]
		Substituting this expression into
		\eqref{osc-rec-upper} and then into \eqref{Hdelta}, we obtain
		\begin{align*}
			H_k(\delta)\leq\sum_{m:C_{km}=1}|\alpha_k|a_k^{s-1}\frac{v_m}{v_k}	H_m(\delta/a_k)	+\frac{C}{v_k}\delta^{s-1}.
		\end{align*}
		Let $C_1=C\max_{1\leq k\leq N}\frac{1}{v_k}.$ Since $\sum_{m=1}^{N}C_{km}|\alpha_k|a_k^{s-1}\frac{v_m}{v_k}=
		\frac{(B(s)v)_k}{v_k} = \theta.$ Therefore,
		\begin{equation}\label{H-rec}
			\begin{split}
				H_k(\delta) \leq\theta\max_{m:C_{km}=1}H_m(\delta/a_k)+C_1\delta^{s-1}.
			\end{split}
		\end{equation}
		Now we will show that $H_k(\delta)$ is uniformly bounded for sufficiently small $\delta$. Let $a_-=\min_{1\leq k\leq N}a_k, \quad a_+=\max_{1\leq k\leq N}a_k.$ 
		Then, $0<a_-\leq a_k \leq a_+<1.$ Let $\delta_0>0$ and define $\delta_1=a_-\delta_0.$ For every $\delta\in[\delta_1,\delta_0]$, the  mesh $\mathcal{P}_{\delta,\tau}(I_k)$ contains at most $\frac{|I_k|}{\delta_1}+2$ nonempty intervals. Moreover, for every $Q \in \mathcal{P}_{\delta,\tau}(I_k)$, $Q \subseteq I$, $ \operatorname{osc}(f,Q)\leq\operatorname{osc}(f,I)$. Therefore,
		\[O_k(\delta)\leq\left(\frac{|I_k|}{\delta_1}+2\right)\operatorname{osc}(f,I).
		\]
		Substituting it into \eqref{Hdelta} yields
		\[H_k(\delta)\leq\frac{\delta_0^{s-1}}{v_k}\left(\frac{|I_k|}{\delta_1}+2\right)\operatorname{osc}(f,I)\leq K_k,\]
		for every $k=1,\dots,N$ and every $\delta\in[\delta_1,\delta_0]$, where $K_k>0$ is a finite constant. Set $K_0=\max_{1\leq k\leq N}K_k$, hence $ H_k(\delta)\leq K_0$ for every $\delta\in [\delta_1,\delta_0]$. Now choose $K\geq\{K_0,\frac{C_1\delta_0^{s-1}}{1-\theta}\}$. Then, $\theta K+C_1\delta_0^{s-1}\leq K$. We claim that
		\begin{equation*}\label{eq:H-bound}
			H_k(\delta)\leq K,
		\end{equation*}
		for every $k=1,\dots,N$ and every $0<\delta\leq\delta_0$. Since the claim holds for $0\leq \delta_1 \leq \delta\leq \delta_0$. We only need to prove this claim for $0<\delta\leq \delta_1$. We prove this using induction, for this suppose that for some integer $n\geq0$, it holds whenever $a_+^n\delta_1\leq\delta\leq\delta_0.$
		Let
		\[a_+^{n+1}\delta_1\leq\delta\leq a_+^n\delta_1.\]
		For every $k$,
		\[\frac{\delta}{a_k}\geq\frac{\delta}{a_+}\geq a_+^n\delta_1.\]
		Since $\delta\leq\delta_1=a_-\delta_0$, gives
		\[\frac{\delta}{a_k}\leq\frac{\delta_1}{a_-}=\delta_0.\]
		Thus, $a_+^n\delta_1\leq\frac{\delta}{a_k}\leq\delta_0.$ By the induction hypothesis,
		\[H_m(\delta/a_k)\leq K,\]
		whenever $C_{km}=1$. Hence, 
		\begin{equation*}
			H_k(\delta)\leq\theta K+C_1\delta^{s-1}\leq\theta K+C_1\delta_0^{s-1}\leq K.
		\end{equation*}
		Since $a_+^n\delta_1\longrightarrow0\quad\text{as }n\longrightarrow\infty,$ yields that
		\[H_k(\delta)\leq K,\]
		for every $k$ and every sufficiently small $\delta>0$. It follows that $O_k(\delta)\leq Kv_k\delta^{1-s}.$
		Since there are only finitely many indices, there exists a constant
		$C_2>0$ such that
		\begin{equation}
			O_k(\delta)\leq C_2\delta^{1-s}, \qquad k=1,\dots,N.
		\end{equation}
		By Lemma \ref{lem:osc-cover}, we have 
		\begin{align*}
			N_\delta(G_k) &\leq A_1\left(\frac{|I_k|}{\delta}+C_2\delta^{-s}+1\right).
		\end{align*}
		Since $s>1$, for all sufficiently small $\delta$, $\delta^{-1}\leq\delta^{-s}.$ Therefore, there exists a constant $C_3>0$, independent of $\delta$, such that
		\[N_\delta(G_k)\leq C_3\delta^{-s}.\]
		Summing over $k$ yields
		\begin{align*}
			N_\delta(G(f))&\leq\sum_{k=1}^{N}N_\delta(G_k) \leq C'\delta^{-s}.
		\end{align*}
		Consequently, $\overdimB G(f)\leq s.$ Since this holds for every $s>D_*$, letting $s\downarrow D_*$ gives $\overdimB G(f)\leq D_*.$
	\end{proof}

	For a word $\omega=(k_1,\ldots,k_n),  k_i\in\{1,\ldots,N\}$, satisfying \(C_{k_i k_{i+1}}=1\) for \(i=1,\dots,n-1\). The total horizontal contraction along the word $\omega$ is $a_\omega=a_{k_1}a_{k_2}\cdots a_{k_n}$ and the total vertical contraction along the word $\omega$ is $\alpha_\omega=\alpha_{k_1}\alpha_{k_2}\cdots\alpha_{k_n}.$ Let $\tau$ and $\omega$ be two words. Then, $\omega$ is said to be less than $\tau$, i.e., $\omega \prec \tau$, if $\tau=\omega \eta$ for some word $\eta$.
	Further, let $W_\omega = W_{k_1}\circ W_{k_2}\circ\cdots\circ W_{k_n}.$
	The corresponding cylinder graph piece is defined by
	\[
	G_\omega=W_\omega(G_{J_{k_n}}).
	\] 
	Define $I_\omega=L_{k_1} \circ \cdots \circ L_{k_n}(J_{k_n})$, where $I_\omega$ is the interval on the x-axis over which the graph $G_\omega$ lies. Therefore, $|I_\omega| = a_\omega |J_{k_n}|$. Consider $J_{\min}= \min_{k}|J_k| >0$ and $J_{\max}= \max_{k}|J_k|$. Therefore, 
	\begin{equation}\label{jminmax}
		J_{\min} a_\omega \leq |I_\omega| \leq J_{\max}a_\omega.
	\end{equation}
	\begin{definition}
		Let $\omega=(k_1,\ldots,k_n) $ and $
		\tau=(\ell_1,\ldots,\ell_m)$ be two words. We say that $\omega$ is an initial subword of $\tau$ if $n\le m$ and $k_i=\ell_i,~ i=1,\ldots,n.$ If neither $\omega$ is an initial subword of $\tau$ nor $\tau$ is an initial subword of $\omega$, then $\omega$ and $\tau$ are said to be incomparable.
	\end{definition}
	\begin{lemma}\label{separationlemma}
		Let $\omega=(k_1,\ldots,k_n)$ and $\tau=(\ell_1,\ldots,\ell_m)$ be two distinct words which are incomparable. Then, the corresponding intervals have disjoint interiors, that is, $I_\omega^\circ\cap I_\tau^\circ=\varnothing.$ In particular, intervals corresponding to distinct words of the same length have disjoint interiors.
	\end{lemma}
	\begin{proof}
		Since $\omega$ and $\tau$ are incomparable, there exists a first
		index $r$ at which they differ. Thus, $k_i=\ell_i,~ i=1,\ldots,r-1,$ and $	k_r\neq \ell_r.$ Since the basic intervals $I_1,\ldots,I_N$ have disjoint interiors, we have
		\[I_{k_r}^{\circ}\cap I_{\ell_r}^{\circ}=\varnothing.	\]
		If $r=1$, then by the nesting property, $I_\omega\subseteq I_{k_1},$ and $I_\tau\subseteq I_{\ell_1}.$ Hence,
		\[I_\omega^\circ\cap I_\tau^\circ=\varnothing.	\]
		Now suppose that $r\geq2$. Since the two words agree up to the
		$(r-1)$th position. Define 
		\[P=L_{k_1}\circ\cdots\circ L_{k_{r-1}}.\]
		Since $I_{k_r}\subseteq J_{k_{r-1}}$ and $	I_{\ell_r}\subseteq J_{k_{r-1}}$, and  each $L_k$ is an affine homeomorphism. Therefore, the map $P$ is one-to-one yields.
		\[P(I_{k_r})^\circ \cap P(I_{\ell_r})^\circ=\varnothing.\]
		Moreover, by the nesting property, $I_\omega\subseteq P(I_{k_r})$ and
		$ I_\tau\subseteq P(I_{\ell_r}).$ Consequently,
		\[I_\omega^\circ\cap I_\tau^\circ=\varnothing.\]
		This completes the proof.
	\end{proof}
	
	We impose the following non-degeneracy condition to obtain the lower bound for the box dimension.
	
	\begin{condition}\label{ass:nondeg}
		There exist constants $c_0>0$ and $r_0>0$ such that for every word $\omega$ satisfying $a_\omega<r_0$,
		\begin{equation}\label{eq:cylinder-lower}
			\operatorname{osc}(f,I_\omega)\geq c_0|\alpha_\omega|,
		\end{equation}
		
	\end{condition}
	
	This condition prevents excessive cancelation of oscillations.

	\begin{lemma}\label{lem:lower-cylinder}
		Let Condition \ref{ass:nondeg} holds. Then, for every
		$1<s<D_*$, there exists a constant $C_s>0$ such that, for all
		sufficiently small $\delta>0$,
		\begin{equation}\label{eq:lower-cover-cylinder}
			N_\delta(G)\geq C_s\delta^{-s}.
		\end{equation}
	\end{lemma}
	
	\begin{proof}
		Let $\mathcal W(\delta)$ be the collection of words satisfying $a_\omega\leq\delta<a_{\omega^-}$, where $\omega^-$ is the word obtained from $\omega$ by deleting its final symbol, i.e., $\omega^- = (k_1,k_2,\cdots,k_{n-1})$ with total horizontal contraction $a_{\omega^-}=a_{k_1}a_{k_2}\cdots a_{k_{n-1}}$ and $a_{\omega^-}=\frac{a_\omega}{a_{k_n}}$.
		Since $a_{\min}=\min_k a_k>0$, every $\omega\in\mathcal W(\delta)$ satisfies
		\begin{equation}\label{eq:stopping-comparable}
			a_{\min}\delta < a_\omega\leq\delta.
		\end{equation}
		We first observe that any two distinct words in $\mathcal W(\delta)$ are incomparable with respect to $\prec$. Suppose that $\omega,\tau\in\mathcal W(\delta)$ and that $\omega\prec\tau$. Since $\omega\in\mathcal W(\delta)$, we have $a_\omega\leq\delta.$ As $\tau$ is a proper extension of $\omega$ and $0<a_k<1$ for every $k$, it follows that $a_{\tau^-}\leq a_\omega\leq\delta.$ However, $\tau\in\mathcal W(\delta)$, which gives $\delta<a_{\tau^-},$ which is a contradiction. Hence, any two distinct words in $\mathcal W(\delta)$ are incomparable. Therefore, by Lemma \ref{separationlemma}, we obtain
		$I_\omega^\circ\cap I_\tau^\circ=\varnothing.$ Moreover, using \eqref{jminmax} and \eqref{eq:stopping-comparable}, we have
		\[a_{\min}J_{\min}\delta<|I_\omega|\leq J_{\max}\delta.\] 
		Thus, the intervals $I_\omega$, $\omega\in\mathcal W(\delta)$ have lengths comparable to $\delta$ and disjoint interiors. Moreover, for a corresponding graph $G_\omega$, Condition \ref{ass:nondeg} gives a vertical height at least $c_0|\alpha_\omega|$. Therefore, square of side $\delta$ covers the vertical distance at most a constant multiple of $\delta$, i.e., 
		\[N_\delta(G_\omega)\geq c\frac{|\alpha_\omega|}{\delta}.\]
		Now, Summing over the $\mathcal W(\delta)$, yields
		\begin{equation*}
			N_\delta(G)
			\geq c\delta^{-1}\sum_{\omega\in\mathcal W(\delta)}|\alpha_\omega|.
		\end{equation*}
		Using \eqref{eq:stopping-comparable}, we have
		\[
		|\alpha_\omega|
		= |\alpha_\omega|a_\omega^{s-1}a_\omega^{1-s}
		\geq C\delta^{1-s}|\alpha_\omega|a_\omega^{s-1}.
		\]
		Hence
		\begin{equation}\label{Ndelta}
			N_\delta(G)
			\geq C' \delta^{-s}
			\sum_{\omega\in\mathcal W(\delta)}|\alpha_\omega|a_\omega^{s-1}.
		\end{equation}
		Now fix $1<s<D_*$. Then, by the definition of $D_*$, we have $	\rho=\rho(B(s))>1.$ Let $q_k=|\alpha_k|a_k^{s-1},~ k=1,\ldots,N.$
		Thus,
		\[B_{km}(s)=C_{km}q_k.\]
		Since $B(s)$ is a nonnegative irreducible matrix, the Perron--Frobenius theorem provides the existence of a positive vector $v=(v_1,\ldots,v_N)^T>0$
		such that $B(s)v=\rho v.$	For a word $\omega=(k_1,\ldots,k_n),$
		define
		\[	Q(\omega)=\prod_{i=1}^{n-1}q_{k_i},	\quad	\text{such that~~~} m(\omega)=Q(\omega)v_{k_n}.	\]
		If $\omega j$ is an extension of $\omega$, then
		$C_{k_nj}=1$ and $Q(\omega j)=Q(\omega)q_{k_n}.$
		Hence,
		\begin{equation*}
			\sum_{j:C_{k_nj}=1}m(\omega j) =Q(\omega)q_{k_n} \sum_{j:C_{k_nj}=1}v_j
			=Q(\omega)(B(s)v)_{k_n}
			=	\rho Q(\omega)v_{k_n}
			=\rho m(\omega).
		\end{equation*}
		Since $\rho>1$, we have  
		\[ \sum_{j:C_{k_nj}=1}m(\omega j)  >   m(\omega).\]
		Since $	a_{\max}=\max_{1\leq k\leq N}a_k<1,$ we have
		\[	a_\omega\leq a_{\max}^{|\omega|}	\longrightarrow 0	\qquad\text{as }|\omega|\longrightarrow\infty.	\]
		Hence, every infinite sequence of word extensions eventually contains a word belonging to $\mathcal W(\delta)$ and $\mathcal W(\delta)$ is finite.
		Starting from the one-letter words and repeatedly replacing each word which does not belong to $\mathcal W(\delta)$ by all of its extensions, we obtain
		\begin{equation}\label{eq:mass-cutset}
			\sum_{\omega\in\mathcal W(\delta)}m(\omega)
			\geq
			\sum_{k=1}^{N}v_k.
		\end{equation}
		As we have, $|\alpha_\omega|a_\omega^{s-1}	=	\prod_{i=1}^{n}q_{k_i}$. Then,
		\[	m(\omega)	=	|\alpha_\omega|a_\omega^{s-1}	\frac{v_{k_n}}{q_{k_n}}.	\]
		Since $v_k>0$ and $q_k>0$ for every $k$ and there are only finitely many indices. Define $C'' = \max_{1 \leq k \leq N}\frac{v_k}{q_k}$ such that
		\[	m(\omega)	\leq	C''|\alpha_\omega|a_\omega^{s-1},	\]
		for every word $\omega$. Therefore,
		\[ \sum_{\omega\in\mathcal W(\delta)} m(\omega) \leq C'' \sum_{\omega\in\mathcal W(\delta)}	|\alpha_\omega|a_\omega^{s-1}.\]
		Hence, by \eqref{eq:mass-cutset}, there exists a constant $c_s>0$ independent of $\delta$ such that
		\begin{equation}\label{eq:cutset-mass-lower}
			\sum_{\omega\in\mathcal W(\delta)}
			|\alpha_\omega|a_\omega^{s-1}
			\geq c_s.
		\end{equation}
		Finally by \eqref{Ndelta}, we have the desired result
		\[N_\delta(G)\geq	C_s\delta^{-s},\]
		for all sufficiently small $\delta>0$, where $C_s>0$ is
		independent of $\delta$. This completes the proof.
	\end{proof}
	
	\begin{theorem}\label{exactboxdimension}
		Let f be an affine PFIF with $0  < a_k <1 $ and conditions \ref{markov} and  \ref{ass:nondeg} are satisfied. Then, $\dimB G(f)=D_*.$ Equivalently,
		\[	\dimB G(f)=
		\begin{cases}
			1,& \rho\bigl((C_{km}|\alpha_k|)\bigr)\leq1,\\[6pt]
			D,& \rho\bigl((C_{km}|\alpha_k|a_k^{D-1})\bigr)=1.
		\end{cases}\]
	\end{theorem}
	
	\begin{proof}
		The upper bound $\overdimB G(f)\leq D_*$ is Theorem \ref{thm:upper-box}. 
		Assume now $D_*>1$. Let $1<s<D_*$. Then, $\rho(B(s))>1$. By Lemma \ref{lem:lower-cylinder}, we have $N_\delta(G)\geq C_s\delta^{-s}$ and for all sufficiently small $\delta$. Therefore
		\[\underdimB G(f)	=\liminf_{\delta \to 0}\frac{\log N_\delta(G)}{-\log\delta} \geq s.\]
		Letting $s\uparrow D_*$ gives $\underdimB G(f)\geq D_*$. Hence together with the upper bound yields $\dimB G(f)=D_*.$
	\end{proof}
	
	\begin{remark}
		The conditionn \ref{ass:nondeg} guarantees that the oscillation of each graph piece remains comparable to the vertical scaling factor $\alpha_\omega$. Without this condition, oscillations may cancel out during the iterative process, leading to an overestimation of the box dimension. Consequently, the upper bound remains valid, but the lower-bound argument and hence the equality $\dimB G(f)=D_*$ may fail.
	\end{remark}
	
	\begin{example}\label{ex}
		Suppose $I=[0,3]$, the interpolation data
		$(t_0,y_0)=(0,0),~(t_1,y_1)=\left(1,\frac{7}{2}\right),~(t_2,y_2)=(2,5)$, and $(t_3,y_3)=\left(3,\frac{5}{2}\right)$. Let $I_1=[0,1], ~I_2=[1,2],~ I_3=[2,3]$ and $J_1=[0,2], ~J_2=[1,3],~J_3=[0,2].$ Then, $J_1=I_1\cup I_2,~J_2=I_2\cup I_3$, and $J_3=I_1\cup I_2$. Therefore, the horizontal contraction ratio are $a_1=a_2=a_3=\frac12$. In particular, $0<a_k<1$ for every $k=1,2,3$. Define the affine maps by
		\begin{equation*}
			L_1(t)=\frac{t}{2} \quad t\in J_1; \qquad
			L_2(t)=\frac{t+1}{2} \quad t\in J_2; \qquad
			L_3(t)=2+\frac{t}{2} \quad t\in J_3
		\end{equation*}
		Consider the vertical maps
		\[M_1(t,y)=-\frac18t+\frac34y; \quad
		M_2(t,y)=\frac12t-\frac12y+\frac{19}{4}; \quad
		M_3(t,y)=\frac58t-\frac34y+5\]
		Thus, $\alpha_1=\frac34$, $\alpha_2=-\frac12$,	$\alpha_3=-\frac34$. The corresponding connection matrix is
		\[C=\begin{pmatrix}
			1&1&0\\
			0&1&1\\
			1&1&0
		\end{pmatrix}.\]
		This matrix is irreducible. Hence, Condition \ref{markov} is satisfied. To verify Condition \ref{ass:nondeg}, define $g(t)=f(t)-\frac12t$ and the interpolation values of $g$ are $g(0)=0$, $g(1)=3$, $g(2)=4$, and $g(3)=1$. Moreover, the corresponding equations for $g$ become
		\[g(L_1(t))=\frac34g(t); \quad 
		g(L_2(t))=-\frac12g(t)+\frac92; \quad
		g(L_3(t))=-\frac34g(t)+4\]
		Thus, the dependence on the horizontal variable disappears. Therefore, iterating along a word $\omega=(k_1,\ldots,k_n)$ gives
		\[ g\bigl(L_{k_1}\circ\cdots\circ L_{k_n}(t)\bigr)= \alpha_\omega g(t)+e_\omega,\]
		where
		$\alpha_\omega=\alpha_{k_1}\alpha_{k_2}\cdots\alpha_{k_n}$
		and $e_\omega$ is a constant. Since, $ I_\omega  =       L_{k_1}\circ\cdots\circ L_{k_n}(J_{k_n}).$ Therefore, 
		\[ \operatorname{osc}(g,I_\omega) = |\alpha_\omega|\operatorname{osc}(g,J_{k_n}). \]
		Hence,
		\begin{align*}
			\operatorname{osc}(g,J_1)&\geq |g(2)-g(0)|=4,\\
			\operatorname{osc}(g,J_2)&\geq |g(2)-g(3)|=3,\\
			\operatorname{osc}(g,J_3)&\geq4.
		\end{align*}
		It follows that $\operatorname{osc}(g,I_\omega)\geq3|\alpha_\omega|$. Since $f(t)=g(t)+\frac12t$, we have $\operatorname{osc}(f,I_\omega)	\geq\operatorname{osc}(g,I_\omega)-\frac12|I_\omega|$.
		Moreover, since $|J_k|=2$ for every $k$,
		$|I_\omega|=a_\omega|J_{k_n}|=2a_\omega$.
		Therefore,
		$\operatorname{osc}(f,I_\omega)
		\geq3|\alpha_\omega|-a_\omega$. Since $a_k=\frac12\leq|\alpha_k|$ for every $k$, it follows that $a_\omega\leq|\alpha_\omega|$. Hence
		\[\operatorname{osc}(f,I_\omega)\geq2|\alpha_\omega|.\]
		Thus, Condition \ref{ass:nondeg} is satisfied with $c_0=2$ and any
		fixed $0<r_0<1$. Since $a_1=a_2=a_3=\frac12$, we obtain
		\[B(s)=2^{1-s} 
		\begin{pmatrix}
			3/4&3/4&0\\
			0&1/2&1/2\\
			3/4&3/4&0
		\end{pmatrix}.\]
		Let $B(s)=2^{1-s} A$. Then, the characteristic polynomial of $A$ is
		$\lambda^2\left(\lambda-\frac54\right)$.
		Hence, $\rho(A)=\frac54$ and therefore	$\rho(B(s))=\frac54\,2^{1-s}$. In particular, $\rho(B(1))=\frac54>1$ and $\rho(B(2))=\frac58<1$.
		Thus, there exists a unique $D_*\in(1,2)$ satisfying
		$\rho(B(D_*))=1$. Therefore, \[D_*=1+\log_2\left(\frac54\right)\approx1.321928.\]
	\end{example}

	\subsection{Classical FIF as a special case}
	
	If $J_k=I$ for every $k$, then $C_{km}=1$ for every $k,m$. In that case, $\sum_{k=1}^{N}|\alpha_k|a_k^{D-1}=1.$
	If the interpolation points are equally spaced, then $a_k=1/N$ and $\sum_{k=1}^{N}|\alpha_k|N^{1-D}=1.$ Therefore, if $\sum_{k=1}^{N}|\alpha_k|>1$, then
	\begin{equation}\label{eq:classical-equal}
		\dimB G(f)=1+\frac{\log\left(\sum_{k=1}^{N}|\alpha_k|\right)}{\log N}.
	\end{equation}
	If $\sum_{k=1}^{N}|\alpha_k|\leq1$, then $\dimB G(f)=1.$
	
	\section{Conclusion and Future Directions}\label{sec8}
	
	In this paper, we have extended the framework of FIFs with partial self-similarity in both univariate and bivariate settings. First, we have constructed a class of bivariate fractal interpolation surfaces on rectangular domains and established the existence of the corresponding fractal interpolation surfaces with partial self-similarity (Lemma \ref{Lemma1}), (Theorem~\ref{Th1}), and (Theorem~\ref{Th2}). The proposed construction generalizes several existing models of fractal interpolation by incorporating local self-similar structures. Furthermore, we have introduced the associated $\alpha$- partial fractal functions. For the univariate case, we have established the existence of the corresponding $\alpha$- partial fractal functions (Theorem~\ref{Th3}) together with their error estimates (Theorem~\ref{Th4}) and extended the construction to the bivariate case (Theorem~\ref{Th5}) and obtained the corresponding error estimates (Theorem~\ref{errorbivariate}). In addition, we investigated the stability and sensitivity of the proposed FIFs with respect to perturbations of the scaling parameters and the underlying interpolation data.
	
	Moreover, we have investigated invariant probability measures associated with the underlying IFSs (Theorem \ref{Th6}) and established the existence and uniqueness of invariant measures supported on the graphs of the corresponding FIFs for the univariate case. Finally, we have studied the box-counting dimension of the graph of the univariate PFIF and obtained upper and lower estimates. Under suitable additional conditions, these estimates lead to an exact formula for the box dimension (Theorem \ref{exactboxdimension}). The resulting dimension formula agrees with the corresponding self-affine RFIF result, while the proof is developed directly within the PFIF framework. As a special case, dimension formula reduces to the well-known box-dimension formula for classical affine FIFs.
	
	The results obtained in this paper open several possible directions for future research. One natural extension is the investigation of the dimension of the proposed bivariate fractal interpolation surfaces. Another interesting problem is the study of regularity properties such as Hölder continuity, differentiability, and bounded variation for partially self-similar fractal interpolation functions and surfaces. It would also be worthwhile to analyze the approximation properties and convergence behavior of the associated $\alpha$ fractal operators in various function spaces.
	
	\section*{Declaration}\noindent
	\textbf{Conflicts of interest.} We do not have any conflict of interest.\\
	\\
	\noindent
	\textbf{Data availability:} No data were used to support this study.\\
	\\
	\noindent
	\textbf{Code availability:} Not applicable\\
	\\
	\noindent
	\textbf{Authors' contributions:} P. Singh: Conceptualization, Writing – original draft, Writing – review \& editing.\\
	M. Pandey: Conceptualization, Writing – original draft, Writing – review \& editing.\\
	S. K. Katiyar: Supervision, Writing – review \& editing.\\
	T. Som: Supervision, Writing – review \& editing.
	
	\section*{Acknowledgements}\noindent
	This work is supported by a Postdoctoral fellowship from Northwest University, Xi'an, Shaanxi Province, China, to the second author. Further, supported by UGC $(231620024614)$ to the first author.
	
	\bibliographystyle{elsarticle-num}
	\bibliography{Reference}
	
\end{document}